\newcommand{\calE}{\mathcal{E}}

\newcommand{\calH}{\mathcal{H}}

\newcommand{\calO}{\mathcal{O}}

\newcommand{\calR}{\mathcal{R}}

\newcommand{\frakm}{\mathfrak{m}}
\newcommand{\frakA}{\mathfrak{A}}
\newcommand{\frakS}{\mathfrak{S}}

\newcommand{\bbA}{\mathbb{A}}
\newcommand{\bbF}{\mathbb{F}}
\newcommand{\bbC}{\mathbb{C}}
\newcommand{\bbH}{\mathbb{H}}

\newcommand{\bbP}{\mathbb{P}}
\newcommand{\bbQ}{\mathbb{Q}}

\newcommand{\bbZ}{\mathbb{Z}}
\newcommand{\bbG}{\mathbb{G}}

\def\SL{{\textrm{SL}}}
\def\Aut{{\text{Aut}}}
\def\AGL{{\text{AGL}}}

\def\Pic{{\text{Pic}}}
\def\Ker{{\text{Ker}}}

\def\PSU{{\textup{PSU}}}

\def\Cl{{\textrm{Cl}}}
\newcommand{\Spec}{\text{Spec}}
\newcommand{\depth}{\text{depth}}
\newcommand{\DR}{\text{DR}}
\newcommand{\et}{\text{et}}
\newcommand{\crys}{\text{crys}}
\newcommand{\disc}{\text{disc}}
\newcommand{\ord}{\text{ord}}

\newcommand{\Br}{\text{Br}}
\newcommand{\Hom}{\text{Hom}}

\newcommand{\GL}{\text{GL}}
\newcommand\bull{\text{\tiny \textbullet}}
\newcommand{\red}{\text{red}}

\documentclass[11pt]{amsart}
\usepackage{amscd, amsmath, amsthm, amssymb}
\newtheorem{theorem}{Theorem}[section]
\newtheorem{lemma}[theorem]{Lemma}
\newtheorem{proposition}[theorem]{Proposition}

\newtheorem{corollary}[theorem]{Corollary}
\theoremstyle{definition}     

\newtheorem{example}[theorem]{Example}

\theoremstyle{remark}
\newtheorem{remark}[theorem]{Remark}
\numberwithin{equation}{section}

\newcommand{\rank}{\text{rank}}

\begin{document}
\title[Symplectic automorphisms]
 {Finite groups of Symplectic automorphisms of K3 surfaces in positive characteristic}
\author[I. Dolgachev]{Igor V. Dolgachev}
\address{Department of Mathematics, University of Michigan, Ann Arbor, MI 48109, USA}
\email{idolga@umich.edu}
\thanks{Research of the first author is partially supported by NSF grant
DMS-0245203}
\author[J. Keum]{JongHae Keum }
\address{School of Mathematics, Korea Institute for Advanced
Study, Seoul 130-722, Korea } \email{jhkeum@kias.re.kr}
\thanks{Research of the second named author is supported by KOSEF grant R01-2003-000-11634-0}
\begin{abstract}
We show that Mukai's classification of finite groups which may act
symplectically on a complex K3 surface extends to positive
characteristic  $p$ under the assumptions that (i) the order of
the group is coprime to $p$ and (ii) either the surface or its
quotient is not birationally isomorphic to a supersingular K3
surface with Artin invariant 1. In the case without the assumption
(ii) we classify all possible new groups which may appear. We
prove that the assumption (i) on the order of the group is always
satisfied if $p > 11$ and if $p=2,3,5,11$ we give examples of K3
surfaces with finite symplectic automorphism groups of order
divisible by $p$ which are not contained in Mukai's list.
\end{abstract}
\maketitle
\section{Introduction}
A remarkable work of S. Mukai \cite{Mukai} gives a classification
of  finite groups which can act on a complex algebraic K3 surface
$X$ leaving invariant  its  holomorphic 2-form (symplectic
automorphism groups). Any such group turns out to be isomorphic to
a subgroup of the Mathieu group $M_{23}$ which has at least 5
orbits in its natural action on a set of 24 elements. A list of
maximal subgroups with this property consists of 11 groups, each
of these can be realized on an explicitly given K3 surface. A
different  proof of Mukai's result was given later by S. Kond\={o} \cite{Ko}
and G. Xiao  \cite{Xiao} classified all possible topological types of
a symplectic action. Neither Mukai's nor Kond\={o}'s proof extends
to the case of  K3 surfaces over
algebraically closed fields of positive characteristic $p$. In
fact there are known examples of surfaces over a field of positive
characteristic whose automorphism group contains a finite
symplectic subgroup which is not realized as a subgroup of
$M_{23}$ (e.g. the Fermat quartic over a field of characteristic
3, or the surface from \cite{DKo} over a field of characteristic
2).

The main tool used in Mukai's proof is the characterization of the
representation of a symplectic group $G\subset \Aut(X)$  on the
24-dimensional cohomology space
$$H^*(X,\bbQ)=
H^0(X,\bbQ)\oplus H^2(X,\bbQ)\oplus H^4(X,\bbQ).$$ Using the
Lefschetz fixed-point formula and the description of possible
finite cyclic subgroups of $\Aut(X)$ and their fixed-point sets
due to Nikulin \cite {Ni1}, allows one to compute the value
$\chi(g)$ of the character of this representation at any element
of finite order $n$. It turns out that $$\chi(g) = \epsilon(n)$$
for some function $\epsilon(n)$ and the same function describes
the character of the 24-permutation representation of $M_{23}$. A
representation of a finite group $G$ in a finite-dimensional
vector space of dimension 24 over a field of characteristic 0 is
called a Mathieu representation if its character is given by
$$\chi(g) = \epsilon(\ord(g)).$$ The obvious fact that $G$ leaves
invariant an ample class, $H^0(X,\bbQ)$, $H^4(X,\bbQ)$,
$H^2(X,\calO_X)$ and $H^0(X,\Omega_X^2)$, shows that $$\dim
H^*(X,\bbQ)^G\ge 5.$$
These properties and the known
classification of finite subgroups of $\SL(2,\bbC)$ which could be
realized as stabilizer subgroups of $G$ in its action on $X$ shows
that the 2-Sylow subgroups of $G$ can be embedded in $M_{23}$. By
clever and non-trivial group theory arguments  Mukai proves that
subgroups of $M_{23}$ with at least 5 orbits are characterized by
the properties that they admit a rational Mathieu representation
$V$ with $\dim V^G \ge 5$ and their 2-Sylow subgroups are
embeddable in $M_{23}$.

The main difficulties in the study of K3 surfaces over algebraically
closed fields of positive characteristic $p$ arise from the
absence of the Torelli Theorem, the absence of a natural unimodular
integral lattice containing the Neron-Severi lattice, the presence of
supersingular K3 surfaces, and the presence of wild automorphisms.

A group of automorphisms is called {\it wild} if it contains a wild
automorphism, an automorphism of order equal to a power of the characteristic
$p$, and {\it tame} otherwise.

In this paper  we first improve
the results of our earlier paper \cite{DK} by showing that a finite
symplectic group  of automorphisms $G$ is always tame if $p >11$.
Next, we show that  Mukai's proof can be extended to
finite tame symplectic groups in any positive characteristic
$p$ unless both the surface and its quotient are birationally
isomorphic to a supersingular K3 surface with Artin invariant equal to
1 (the exceptional case). To do this, we first
prove that  Nikulin's
classification of finite order elements and its sets of fixed
points extends to positive characteristic $p$, as long as the order is
coprime to $p$. Next we consider the
24-dimensional representations of $G$ on the $l$-adic cohomology,
$l\ne p$,
$$H_{\et}^*(X,\bbQ_l)= H_{\et}^0(X,\bbQ_l)\oplus H_{\et}^2(X,\bbQ_l)\oplus
H_{\et}^4(X,\bbQ_l)$$ and on the crystalline cohomology
$$H_{\crys}^*(X/W)=H_{\crys}^0(X/W)\oplus H_{\crys}^2(X/W)\oplus
H_{\crys}^4(X/W).$$ It is known that the characteristic polynomial
of any automorphism $g$ has integer coefficients which do not
depend on the choice of the cohomology theory. Comparing
$H_{\crys}^2(X/W)$ with the algebraic De Rham cohomology
$H_{\text{DR}}^2(X)$ allows one to find a  free submodule of
$H_{\crys}^*(X/W)^G$ of rank 5 except when both $X$ and $X/G$ are  birationally isomorphic to a supersingular K3 surface
with Artin invariant equal to 1 (the exceptional case). This shows
that in a non-exceptional case, for all prime $l\ne p$, the vector spaces $$V_l =
H_{\et}^*(X,\bbQ_l)$$ are Mathieu representations of $G$ with
$$\dim V_l^G \ge 5.$$ A careful analysis of Mukai's proof shows
that this is enough to extend his proof.

  In the exceptional case, it is known that a supersingular K3
  surface with Artin invariant equal to 1 is unique up to isomorphism.  It is isomorphic to the Kummer
surface of the product of two supersingular elliptic curves if
$p>2$, and to the surface from \cite{DKo} if $p=2$. We call a tame
group $G$ exceptional if it acts on such a surface with $\dim
V_l^G = 4$. We use arguments from Mukai and some additional
geometric arguments to classify all exceptional groups. All
exceptional groups turn out to be subgroups of the Mathieu group
$M_{23}$ with 4 orbits. S. Kond\=o has confirmed the converse,
that is, any subgroup of $M_{23}$ with 4 orbits is either from our
list or wild. The problem of realizing exceptional groups will be
discussed in other publication.

In the last section we give examples  of K3 surfaces in
characteristic $$p = 2,3,5,11$$ with wild finite symplectic
automorphism groups  which are not contained in Mukai's list. We
do not know a similar example in characteristic $p =7$, however we
exhibit a K3-surface with a symplectic group of automorphisms of
order 168 which does not lift to a surface from Mukai's list.

\bigskip
{\it Acknowledgment}
\medskip

We thank S. Kond\={o} for useful discussions and sharing with us
some of the examples of K3 surfaces over a field of small
characteristic with exceptional finite groups of symplectic
automorphisms. We would also like to thank  L. Illusie, T.
Katsura, N. Katz, W. Messing and A. Ogus for an e-mail exchange
and discussion helping us to understand some of the properties of
crystalline and De Rham cohomology. We are grateful to D. Allcock who helped us with a computer computation allowing us to greatly simplify the classification of exceptional groups. Finally we thank
R. Griess for his generous help with explaining to us some group theory.

\section{Automorphisms of order equal to the characteristic $p$}
\bigskip

Let $X$ be a $K3$ surface over an algebraically closed field $k$
of positive characteristic $p$.

The automorphism group Aut($X$) acts on the 1-dimensional space of
regular 2-forms $H^0(X,\Omega^2_X)$. Let
$$\chi_{2,0} : \Aut(X)\to
k^*$$ be the corresponding character. An automorphism $g$ is
called \emph{symplectic} if $\chi_{2,0}(g)=1$.

Obviously any automorphism of order equal to the characteristic is
symplectic. The main result of this section is the following.

\begin{theorem}\label{main} Let $g$ be an  automorphism of $X$ of
  order $p={\rm char}(k)$. Then
$p\le 11$.
\end{theorem}

First, we recall the following result from our previous paper
\cite{DK}.

\begin{theorem} \label{theorem}  Let $g$ be an automorphism of order $p$ and let $X^g$ be the
set of fixed points of $g$. Then one of the following cases occurs.
\begin{itemize}
\item[(1)] $X^g$ is finite and consists of $0, 1$ or $2$ point; $X^g$ may
be empty only if $p=2$, and may consist of $2$ points only if
$p\le 5$.
\item[(2)] $X^g$ is a divisor such that the Kodaira dimension
of the pair $(X,X^g)$ is equal to $0$. In this case $X^g$ is a
connected nodal cycle, i.e. a connected union of smooth rational
curves.
\item[(3)] $X^g$ is a divisor such that the Kodaira dimension
of the pair $(X,X^g)$ is equal to $1$. In this case $p\le 11$ and
there exists a divisor $D$ with support $X^g$ such that the linear
system $|D|$  defines an elliptic or quasi-elliptic fibration
$\phi : X \to \bbP^1$.
\item[(4)] $X^g$ is a divisor such that the Kodaira dimension
of the pair $(X,X^g)$ is equal to $2$. In this case $X^g$ is equal
to the support of some nef and big divisor $D$. Take $D$ minimal
with this property. Let $$d: = \dim H^0(X,{\mathcal
O}_X(D-X^g))\quad {\rm and}\quad N: = {1\over 2}D^2+1.$$ Then
$$p(N-d-1)\le 2N-2.$$
\end{itemize}
\end{theorem}

Theorem \ref{theorem} does not give the bound for $p$ in the cases
(1), (2) and (4). First we take care of the case (4). Its proof is
essentially contained in \cite{DK}.

\begin{proposition} \label{kappa2} Let $g$ be an automorphism of order
$p={\rm char}(k)$ and $F = X^g$ (with reduced structure). Assume
$\kappa(X,F) = 2$. Then $p \le 5$. Moreover, if $p = 5$, then
$X^g$ contains a curve $C$ of arithmetic genus 2 and the linear
system $|C|$ defines a double cover $X \to \bbP^2$. The surface is
$g$-isomorphic to the affine surface
$$z^2 = (y^5-yx^4)P_1(x)+P_6(x), \ g(x,y,z) = (x,y+x,z),$$
where $P_i(x)$ is a polynomial of degree $i$.
\end{proposition}

\begin{proof} We improve the argument from our paper \cite{DK}.

\medskip\noindent
Case 1: $F$ is nef and $|F|$ has non-empty fixed part.

By a well-known result due to Saint-Donat \cite{SD}, $F \sim
aE+\Gamma$, where $a\ge 2$ and $E$ is an irreducible curve of
arithmetic genus 1 and $\Gamma$ is the fixed part which is a
$(-2)$-curve with $E\cdot \Gamma = 1$. In this case $g$ leaves
invariant a genus 1 pencil defined by the linear system $|E|$.
Note that $\Gamma$ is a section of the fibration and is fixed
pointwisely by $g$. It is well known that no elliptic curve admits
an automorphism of order $\ge 5$ fixing the origin. Thus $p\le 3$.

\medskip\noindent
Case 2: $F$ is nef and $|F|$ has no fixed part.

Assume $F^2=0$. Then $F \sim aE$, where $E$ is an irreducible
curve of arithmetic genus 1, hence $\kappa(X,F)=1$.

Assume $F^2=2$. Then $|F|$ has no base points, and $|F|$ defines a
map of degree 2 onto $\bbP^2$. Since we can assume that $p \ne 2$,
the map is separable and the branch curve is a curve of degree 6
with simple singularities which is invariant under a projective
transformation of order $p$. Since $F$ is the pre-image of a line,
$g$ fixes pointwisely a line in $\bbP^2$. Thus it is conjugate to
a transformation $(x_0,x_1,x_2)\mapsto (x_0, x_1,x_2+x_1)$. It is
known that the ring of invariants $k[x_o,x_1,x_2]^g$ is generated
by the three polynomials $x_0,x_1, x_2^p-x_2x_1^{p-1}$. Thus $p\le
5$.

Assume $F^2\ge 4$. Then $F$ is nef and big, hence Theorem
\ref{theorem} (4) applied to $D = F$ gives $d = 1$ and $p(N-2)\le
2N-2$. Since $N = D^2/2+1\ge 3$, this gives $p \le 3$.

\medskip\noindent
Case 3: $F$ is not nef.

Since $F$ is reduced and connected, being non-nef means that
$$F = F_1+C_1+\ldots +C_k,$$
where $C_i$'s are chains of $(-2)$-curves with no common
components and $C_i\cdot F_1 = 1$ and $F_1$ is nef.  In fact, let
$E$ be a $(-2)$-curve such that $F\cdot E < 0$. Then $F = E+F'$,
and $E\cdot F = -2+E\cdot F' > -2$. Hence $E\cdot F' = 1$. This
implies that $E$ is an end-component of $F$ and $F'$ is connected.
If $F'$ is not nef, repeat the same process. Continuing in this
way, we prove the claim.

If $F_1^2 = 0$, then $g$ fixes  $F_1$ and $C_1$, hence fixes
pointwisely a fibre of an elliptic fibration and a section. We
argue as in Case 1 to get $p\le 3$. If $F_1^2 = 2$, we use Case 2
and get $p\le 5$. So let us assume that $F_1^2\ge 4$, so $F_1$ is
big and nef.

Take $D =F+F_1=2F_1+C_1+\ldots +C_k$. It is nef, and indeed
minimal nef. Then
$$N =h^0(D)-1 = {1\over 2}D^2+1 ={1\over 2}(4F_1^2+4k-2k)+1 = 2F_1^2+k+1,$$
$$d = h^0(D-F) = h^0(F_1) ={1\over 2}F_1^2+2 \ge 4.$$
Hence $N = 4d+k-7$ and the inequality of Theorem \ref{theorem} (4)
gives
$$p(N-d-1) = p(3d+k-8) \le 2(N-1) = 8d+2k-16.$$
Hence
$$p \le {(8d+2k-16)\over(3d+k-8)} = 2+{2d\over 3d+k-8}<4.$$
\end{proof}

It remains to consider the cases (1) and (2) of Theorem
\ref{theorem}, that is, the cases where $X^g$ is either one point
or a connected nodal cycle. To do this, we first recall the
following information from \cite{DK}(Lemma 2.1 and Theorem 2.4).

\begin{proposition} \label{dico} The following is true.
\begin{itemize}
\item [(1)] If $X^g$ consists of a point, then $X/(g)$ is either a
rational surface
 with trivial canonical divisor and one isolated elliptic
Gorenstein singularity or a K3 surface with one rational double
point.  The latter case occurs only if $p\le 5$.
\item [(2)] If $X^g$ is a nodal cycle, then $X'/(g)$ is a rational surface
 with trivial canonical divisor with one isolated elliptic
Gorenstein singularity, where $X'$ is the surface obtained from
$X$ by blowing down the nodal cycle $X^g$.
\end{itemize}
\end{proposition}

In the following proposition, we prove that the case (2) of
Theorem \ref{theorem} occurs only if $p=2$.

\begin{proposition}\label{newprop} Suppose $X^{g}$ is a nodal cycle. Then $p = 2$.
\end{proposition}

\begin{proof}
The quotient $Z = X/(g)$ is known to be a rational surface with at
most rational singularities and $-K_Z = (p-1)B$, where $B_{\red}$
is the image of the nodal cycle in $Z$ (see \cite {DK}, section
3). Let $\pi:Z'\to Z$ be a minimal resolution of singularities.
Then $-K_{Z'} = -\pi^*(K_Z)+\Delta' = (p-1)\pi^*(B)+\Delta'$,
where $ \Delta'$ is a positive divisor supported on the
exceptional locus. Obviously $Z'$ is a resolution of singularities
of the surface $X'/(g) $, where $X'$ is obtained from $X$ by
blowing down the nodal cycle. Let $a:Z'\to V$ be the blowing down
to a minimal resolution of $X'/(g) $. Then $-K_V =
a_*((p-1)\pi^*(B))+\Delta')$ is the fundamental cycle of $V$. By
Corollary 3.6 of \cite{DK}, we have $H^1((p-1)B, \calO_{(p-1)B})
\cong k$. Since $Z$ has only rational singularities this implies
that $H^1((p-1)\pi^*(B),\calO_{(p-1)\pi^*(B)}) \cong k,$ and hence
$H^1(a_*(p-1)B,\calO_{a_*(p-1)B}) \ne \{0\}$. It is known that for
any proper part $A$ of the fundamental cycle of a minimal elliptic
singularity, we have $H^1(A,\calO_A) = 0$. This implies that
$$-K_{V} = (p-1)a_*(\pi^*(B)).$$ Since $V$ is a non-minimal
rational surface, it contains a $(-1)$-curve $E$. Intersecting
both sides of the previous equality with $E$, we get $1 =
(p-1)\pi^*(B)\cdot E$, hence $p = 2$.
\end{proof}

\begin{lemma} (P. Samuel)\label{samuel} Let $A$ be a normal noetherian local
$k$-algebra of dimension $\ge 2$ with maximal ideal $\frakm$. Let
$G$ be a finite group of  automorphisms of $A$  with local ring of
invariants $A^G$. Assume that $G$  acts freely on the punctured
local scheme $V = \Spec(A)\setminus \{\frakm\}$. Then the class
group ${\rm Cl}(A^G)$ fits in the following exact sequence of
groups:
$$0\to H^1(G,A^*)\to {\rm Cl}(A^G)\to H^0(G,H^1(V,\calO_V^*)).$$
In particular, if in addition $A$ is factorial,  ${\rm
Cl}(A^G)\cong H^1(G,A^*)$.
\end{lemma}

\begin{proof} Let $U = \Spec(A^G)\setminus \{\frakm^G\}$. It is known that
$\depth(A^G) \ge 2$ \cite{Fogarty}. Thus the class group
$\Cl(A^G)$ is isomorphic to the Picard group $\Pic(U)$. We apply
the two spectral sequences which we used in \cite{DK} to the free
action of $G$ on $V$ with quotient $U$ and the $G$-linearized
sheaf $\calO_V^*$.
$$E_2^{i,j} = H^i(G,H^j(V,\calO_V^*)) \Rightarrow \bbH^n,$$
$${}'E_2^{i,j} = H^i(U,\calH^j(G,\calO_V^*))\Rightarrow \bbH^n.$$
Since $G$ acts freely, $$\calH^j(G,\calO_V^*) = 0\quad {\rm
for}\quad j>0$$ and the second spectral sequence gives an
isomorphism $$\Cl(A^G)= \Pic(U) = H^1(U,\calO_U^*) \cong \bbH^1.$$
Now the first assertion follows from the first spectral sequence.

If $A$ is factorial, then $$\Cl(A) =H^1(V,\calO_V^*) = 0,$$ which
proves the last assertion.
\end{proof}

The well-known property of cohomology groups (\cite{Cartan-Eilenberg},
Chapter XII, Proposition 2.5) gives the following.

\begin{corollary}\label{killed} In the situation of the above lemma, if in addition $|{\rm Cl}(A)|$ is finite, then ${\rm Cl}(A^G)$ is
killed by multiplication by the product $|{\rm Cl}(A)|\cdot |G|$.
\end{corollary}

We will also use the following formula from \cite{Kato}, Lemma
4.1.7 or \cite{Saito}, Theorem 7.4.

\begin{lemma}\label{kato} Let $X$ be a smooth surface over an algebraically closed field $k$ of characteristic $p\ge 0$ and $x\in X$ be a closed point of $X$. Let $G$ be a finite group of automorphisms of $X$ such that $X^G = \{x\}$. Let $U = X\setminus \{x\}$ and $V = U/G$. Then
\begin{equation}\label{katosaito}
e_c(U) = (\#G-1)(e_c(V)+1)+e_c(V)-\sum_{g\in G\setminus\{1\}}l(g),
\end{equation}
where $e_c(Z)$ denotes the l-adic Euler-Poincar\'e characteristic with compact support for any $l\ne p$ and
$l(g)$ is the intersection index of the graph of $g$ with the diagonal at the point $(x,x)$.
\end{lemma}

\begin{proposition}\label{prop1} Let $g$ be a wild  automorphism of order
$p \ne 11$  \footnote{The assumption $p\ne 11$ can be dropped, see \cite{DK2}.} acting on a K3 surface such that $|X^g| = 1$.
Assume that the
quotient surface is rational. Then there exists a $g$-invariant
elliptic or quasi-elliptic fibration on $X$.
\end{proposition}

\begin{proof}
Note that the quotient surface $X/(g)$
is a rational surface with trivial canonical divisor with one
isolated elliptic Gorenstein singularity $Q$. Let
$$\sigma: Y \to X/(g)$$
be a minimal resolution. We have
$$K_{Y}=-\Delta ,$$
where $\Delta$ is an effective divisor whose support is equal to
the exceptional set of $\sigma$ and satisfies
$$\Delta\cdot R_i \le 0$$
for any irreducible component $R_i$ of $\Delta$ (see, for example,
\cite{Reid}, 4.21). Let $N$ be the sublattice of $\Pic(Y)$ spanned
by the divisor classes of the curves $R_i$ and $N^*$ be its dual
lattice. It is known that the class group of the local ring of the
singular point $Q$ is mapped surjectively onto $N^*/N$
(\cite{Gr}). By Corollary \ref{killed}, the group $N^*/N$ is a
$p$-elementary finite abelian group. Since $Y$ is a rational
surface, the Picard lattice $\Pic(Y)$ is a unimodular hyperbolic
lattice, hence the orthogonal complement $N^\perp$ is a
$p$-elementary lattice of signature $(1,t)$, $t\ge 0$. Since $N$
contains $K_Y$, it is also an even lattice. In particular, $\rank
N^\perp\ge 2$ if $p$ is odd.

\medskip\noindent
Case 1: $p>2$ and $\rank N^\perp\ge 3$, or $p=2$ and $\rank
N^\perp\ge 2$.

It follows from \cite{RS} that, if $p> 2$, an even  $p$-elementary
hyperbolic lattice of rank $\ge 3$ is determined uniquely by its
rank and discriminant. An explicit construction of such a lattice
shows that it always contains an isotropic vector. Also it is
known that an even  $2$-elementary hyperbolic lattice of rank $\ge
2$ always contains an isotropic vector (\cite{Ni2}, Theorem
4.3.3). Thus $N^\perp$ contains an isotropic vector $v$. Since the
homomorphism $\sigma^*:\Pic(X/(g)) \to N^\perp$ has finite
cokernel, some multiple of $v$ is equal to $\sigma^*(D)$, where
$D$ is a Cartier divisor class on $X/(g)$ with $D^2 = 0$ (here we
use the intersection theory for $\bbQ$-divisors on
$\bbQ$-factorial surfaces). The pre-image of $D$ in $\Pic(X)$ is
an isotropic vector. Dividing it  by an integer, we can represent
this vector by an effective primitive divisor class
\begin{equation}\label{divA}
A = F+E,
\end{equation}
 where $F$ is the fixed part of $|A|$ consisting of
a bunch of $(-2)$-curves and $|E|$ is a free pencil of arithmetic
genus 1 curves.  Since our automorphism $g$ preserves $|A|$, it
must preserve $|E|$. The arithmetic genus 1 pencil $|E|$ is
$g$-invariant and one of its fibres contains the point $X^g$. This
proves the assertion.

\medskip\noindent
Case 2: $p>2$ and $\rank N^\perp = 2$.

We apply Lemma \ref{kato} to our situation. Since $e_c(X\setminus
\{x\})=23$, the formula \eqref{katosaito} gives
\begin{equation}\label{ks2}
23 = (p-1)(e_c(V) + 1-l(g))+e_c(V).
\end{equation}
Keeping the notation of the lemma we have $V \cong
Y\setminus \Delta_{\text{red}}$. By the additivity of $e_c$ we
have
$$e_c(V) = e_c(Y)-e_c( \Delta_{\text{red}}) = 2+\rank~\Pic(Y)-(1+\rank~N) = \rank~N^\perp +1 = 3.$$
Here we use that the graph of components of $\Delta$ is a tree of
smooth rational curves or an irreducible cuspidal curve of arithmetic genus 1. Otherwise, the local Picard group of
$X/(g)$ contains a connected algebraic group isomorphic to an
elliptic curve or the multiplicative group, hence
contains elements of finite order not killed by multiplication by
$p$. The latter contradicts Corollary \ref{killed}. Now formula
\eqref{ks2} gives $20 = (p-1)a$, where $a < 4$. The only
possibility is  $p = 11$.

\medskip\noindent
Case 3: $p=2$ and $\rank N^\perp = 1$.

In this case $e_c(V)=2$, hence the formula \eqref{ks2} gives
$$21 =(2-1)(2+1-l(g)),$$
absurd.
\end{proof}

\begin{remark} If $\rank~ N^\perp > 2$, the formula \eqref{ks2} gives
$$23-k = (p-1)(k+1-l(g)),$$
where $k = e_c(V)>3$. If $k \ne  23$, this gives a weaker bound $p
\le 19$ with possible cases $(p,k,l(g))=(19,5,5), (17,7,7),
(13,11,11)$, etc.
\end{remark}

\begin{proposition} \label{prop2} Let $g$ be an  automorphism of $X$ of
  order $p={\rm char}(k)$. If $X$ admits a $g$-invariant arithmetic genus
$1$ fibration, then $p\le 11$.
\end{proposition}

\begin{proof} Since quasi-elliptic fibrations occur only in characteristic
$p =2$ or $3$, we may assume that  our fibration is an
elliptic fibration. Assume first that $g^*$ acts as identity on
the base curve $\bbP^1$ of the elliptic fibration. Then $g$
becomes an automorphism of the elliptic curve $X/\bbP^1$ over
the function field of $\bbP^1$. On the Jacobian of this
elliptic curve $g$ induces an automorphism $g'$ of order $p$. Let
$$j:J\to\bbP^1$$ be the jacobian elliptic fibration. Note that $J$ is a $K3$
surface (cf. \cite{CD}, Theorem 5.7.2). Since the order of $g'$ is
$p>3$, $g'$ must be a translation
by a $p$-torsion section. By Corollary 5.9 from \cite{DK},  $p\le
11$.

Next assume that $g^*$ acts non-identically on the base ${\bf
P}^1$ of the elliptic fibration. Let $s_0$ be the unique fixed
point of $g^*$ on the base. The fibre $X_{s_0}$ contains $X^g$.
The remaining singular fibres form orbits of fibres of the same
type. By the same argument as in the proof of Corollary 5.6
\cite{DK}, we see that $p\le 11$.
\end{proof}

Now, Propositions \ref{kappa2}, \ref{newprop}, \ref{prop1} and
\ref{prop2} together with Theorem \ref{theorem} prove our Theorem
\ref{main}.

\begin{remark} One should compare the result of Theorem \ref{main}
  with the known  result that  an abelian surface $A$ over a field of
  characteristic $p > 0$ does not admit an automorphism $g$ of order
  $p > 5$ which fixes a point. This result can be proved by
  considering the action of $g$ on the Tate module
  $H_{\et}^1(A,\bbQ_l), l\ne p$ and applying the Weyl theorem. A
  similar proof for
 K3 surfaces only shows that $p > 23$ is impossible. We thank Yuri Zarkhin for this remark.
\end{remark}

\medskip
If $g$ preserves an elliptic or a quasi-elliptic pencil, then it
either acts non-identically on the base or is realized by an automorphism
 of its
jacobian fibration. In the latter case, if $p>3$, then it is realized
by a translation by a $p$-torsion section.
In \cite{DK} (Corollary 5.9) we have proved
that no non-trivial $p$-torsion section exists if $p> 11$. One can
improve this bound to 7 by the same proof. In the case $p =11$, if an
$11$-torsion section exists, then the proof (see information
(i)-(iii) in \cite{DK}, p. 126) gives only one possible
combination of types of singular fibres: three singular fibres of
type $I_{11}$, $I_{11}$, $II$. If this happens, the formula (5.1)
from \cite{DK} gives a contradiction; in this case the left hand
side of the formula cannot be an integer.

We state this result for future references.

\begin{theorem}\label{corr} Let $f:X\to \bbP^1$ be an elliptic
fibration with a section on a K3-surface over an algebraically
closed field of characteristic $p > 7$.  Then the group of
$p$-torsion sections is trivial.
\end{theorem}

A different proof of the above result is given by A. Schweizer
\cite{Sch}.

\section{Tame symplectic automorphisms}

\bigskip
In this section we consider symplectic automorphisms of $X$ of
finite order  prime to the characteristic $p$. We will show that they
behave as in the complex case.

\begin{lemma}\label{klein} Let $\Gamma$ be a finite subgroup of
  $\SL(2,k)$ of order prime to $p$. Then $\Gamma$ is isomorphic to a finite
  subgroup of $\SL(2,\bbC)$, i.e. one of the following groups:
a cyclic, a binary dihedral(=quaternion), binary tetrahedral, binary
  octahedral,
 or binary icosahedral group.
\end{lemma}

\begin{proof} This is of course well-known. For completeness sake let us recall  the usual proof
(going back to Felix Klein).  Let  $\Gamma'\subset {\rm PSL}(2,k)$ be the image of $\Gamma$ in
${\rm PSL}(2,k)$. Any non-trivial element $g\in \Gamma' $ has exactly 2 fixed
points in the natural action of on $S = \bbP^1(k)$. Let
$\mathfrak{P}$ be the union of the sets of fixed points $S^g, g\in
\Gamma'\setminus \{1\}$. Let $O_1,\ldots,O_r$ be the orbits of
$\Gamma'$ in $\mathfrak{P}$ and $n_1,\ldots, n_r$ be the orders of
the corresponding stabilizers. An easy argument using the Burnside
counting formula, gives the equation
$$\sum _{i=1}^{r}(1-\frac{1}{n_i}) = 2-\frac{2}{|\Gamma'|}.$$
This immediately implies that either $r = 2$ and $\Gamma'$ is a cyclic group, or $r = 3$, and
$$(n_1,n_2,n_3;|\Gamma'|) = (2,2,n;2n), (2,3,3;12), (2,3,4;24),\,\,{\rm or}\,\,(2,3,5;60).$$
An easy exercise in group theory shows that the groups $\Gamma$ are isomorphic to a cyclic or a binary polyhedral group.
\end{proof}

For a nondegenerate lattice $L$, we denote by $\disc(L)$ the
discriminant of $L$. We define
$$d_L:=|\disc(L)|,$$ the order of the discriminant group $L^*/L$.

\medskip
Let $\phi:X\to C$ be an elliptic surface, with or without a section.
For any reducible fibre $X_c$, let  $S_c$ be the sublattice of the Picard lattice $S_X$
generated by all irreducible components of the fibre. The Gram matrix with respect to the basis formed by the irreducible components is described by a Dynkin diagram of affine type $\tilde{A}_n, \tilde{D}_{n}, \tilde{E}_n$, where $n+1$ is the number of irreducible components. The radical of $S_c$
is spanned by the (scheme-theoretical) fibre $X_c$ considered as a divisor on $X$. The
quotient $\bar{S}_c$ by the radical is isomorphic to the corresponding negative definite root lattice of types $A_n, D_{n}, E_n$ (we say that the fibre is of the corresponding type). If $M_c\subset S_c$ is a
negative definite sublattice of maximal rank, then the composition with the projection to $\bar{S}_c$ defines an embedding of lattices
$$M_c\hookrightarrow \bar{S}_c.$$

The orthogonal sum  $\oplus_{c\in C}\bar{S}_c$ is the quotient of the
sublattice  $S_X^{\textup{vert}}$
of $S_X$ generated by components of fibres by the rank 1 sublattice  generated by the divisor class of any fibre. We denote  the orthogonal sum by $\calR(\phi)$ and call it the \emph{root lattice} of the elliptic surface $\phi$.

A negative definite sublattice $M=\oplus_{c\in C}M_c\subset S_X^{\textup{vert}}$ is called \emph{maximal} if  each $M_c$ is a
 sublattice of $S_c$ of maximal possible rank. Its image in $\calR(\phi)$ is a sublattice of finite index, say $a$.  In particular, we have
\begin{equation}\label{max}
d_{M}=a^2 d_{\calR(\phi)}.
\end{equation}

\begin{lemma}\label{square} Let $G$ be a finite group of symplectic
  automorphisms of a K3 surface of order prime to $p$. Let $Y\to X/G$
  be a minimal resolution of the quotient $X/G$ and let $\calR_G$ be
  the sublattice of $\Pic(Y)$ generated by the irreducible components
  of the exceptional divisor. Then the discriminant  of $\calR_G$ is
  coprime to $p$. Moreover, if  $\rank~ \calR_G = 20$, then the discriminant
is not a square.
\end{lemma}

\begin{proof} The lattice $\calR_G$ is a direct sum of irreducible
  root lattices $\calR_i$  of type $A, D, E$ generated by irreducible
  components of a minimal resolution of quotient singularities
  corresponding to stabilizer subgroups $G_i$  of $G$. Via the action
  of $G_i$ on the tangent space of $X$ at one of its fixed points the group $G_i$ becomes
  isomorphic to a finite subgroup $H$ of $\SL(2,k)$. Since $\#G_i$ is
  prime to $p$, the quotient singularity is formally isomorphic to the
  singularity $\bbA_k^2/H$. Now we apply Lemma \ref{klein} and use the
  well-known resolution of the quotient
singularity $\bbA^2/H$.  If $\calR_i$ is of type $A_k$, then $\#G_i = d_{\calR_i} = k+1$. If
$\calR_i$ is of type $D_n$, then $\#G_i = 4(n-2), d_{\calR_i} = 4$. If
$\calR_i$ is of type $E_6$, then $\#G_i = 24, d_{\calR_i} = 3$. If
$\calR_i$ is of type $E_7$, then $\#G_i = 48, d_{\calR_i} = 2$. If
$\calR_i$ is of type $E_8$, then $\#G_i = 120, d_{\calR_i} = 1$. In
all cases we see that if $p|d_{\calR_i}$, then $p|\#G_i$, a
  contradiction to
the assumption that $\#G$ is prime to $p$.
This proves the first assertion.

Assume that $\rank~\calR_G = 20$. Obviously  the Picard number
$\rho$ of $Y$ satisfies $\rho\ge 21$. Thus $Y$ is a supersingular K3
surface. It is known (\cite{Artin3}) that the
discriminant group of the Picard lattice  of a supersingular K3
surface is a $p$-elementary abelian group $(\bbZ/p)^{2\sigma}$,
where $\sigma$ is the Artin invariant of the surface. Let $N$ be
the orthogonal complement of $\calR_G$ in $S_{Y}$. We have
$$d_{\calR_G}\cdot d_N =
i^2p^{2\sigma},$$ where $i$ is the index of the sublattice
$$\calR_G\oplus N\subset S _{Y}.$$

Assume that $d_{\calR_G}$ is a square. Then  $d_N$ is a square. The lattice
$N$ is an indefinite lattice of rank 2 whose discriminant is the
negative of a
square. It must contain a primitive isotropic vector. By
Riemann-Roch, we can represent it by an effective divisor $A$ with
self-intersection 0. Write $A$ as in \eqref{divA}. It is known that
a suitable composition of reflections with respect to the divisor
classes  of $(-2)$-curves sends $A$ to $E$.
Let
$$\psi:S _{Y}\to S _{Y},\,\,\, \psi(A)=E$$
be the composition. Let $R_i, i = 1,\ldots, 20,$ be the irreducible
components of the exceptional divisor. Since $\psi$ is an isometry,
the images
$\psi(R_i)$ generate a sublattice $M$ of $S_Y$ isomorphic to $\calR_G$.
Since $A\cdot R_i=0$ for all $i$, we have
$$E\cdot \psi(R_i)=0,\ i = 1,\ldots,20.$$
Since, by Riemann-Roch, each $\psi(R_i)$ is effective or anti-effective,
this implies that
$$E\cdot C_{ij}=0$$
for all irreducible components $C_{ij}$ of $\psi(R_i)$ or
$-\psi(R_i)$, that is,  $C_{ij}$'s are irreducible components of
divisors from the  pencil  $|E|$ of genus 1 curves.
Since $M$ is a negative definite lattice of rank 20, it is a maximal
sublattice of $S_X^{\textup{vert}}$ with respect to the elliptic
fibration $\phi$ given by the  pencil.  By \eqref{max}
\begin{equation}\label{max2}
a^2d_{\calR(\phi)}=d_{M}=d_{\calR_G}.
\end{equation}
Let
$$f:J\to \bbP^1$$ be the jacobian fibration of $\phi$. The surface
$J$ is a K3 surface with the same type of singular fibres
(\cite{CD}, Theorem 5.3.1). In particular,
$$\calR(\phi)\cong \calR(f).$$
The surface $J$ is a
supersingular K3 surface, since $S_J^{\textup{vert}}$ is of rank 21.
The orthogonal complement of the root lattice $\calR(f)$ in $S_J$ is generated by
the divisor classes of the zero-secton and a fibre, and hence
unimodular. Thus we obtain
$$d_{S_J}m^2 = d_{\calR(f)}$$
for some number $m$ (equal to the order of the Mordell-Weil group of sections of $f$).
Since $d_{S_J} =p^{2\sigma'}$, we obtain that
$p|d_{\calR(f)}$, and, by \eqref{max2}, $p|d_{\calR_G}$. This contradicts the first assertion.
\end{proof}

\begin{theorem} \label{theorem2} Let $g$ be a symplectic automorphism of finite order $n$.
If $(n,p)=1$, then $g$ has only finitely many fixed points $f$ and the
possible pairs $(n,f)$ are as follows.
$$(n,f)=(2,8),(3,6),(4,4),(5,4),(6,2),(7,3),(8,2).$$
\end{theorem}

\begin{proof} At a fixed point of $g$, $g$ is linearizable because
$(n,p)=1$. This implies that the quotient surface $Y=X/(g)$ has at
worst cyclic quotient Gorenstein singularities  and its minimal
resolution is a $K3$ surface. So, Nikulin's argument \cite{Ni1} for
the complex case works (when  we replace the rational cohomology
with the $l$-adic cohomology) except in the following two cases.

Case 1: $n=11$ and $X/(g)$ has two $A_{10}$-singularities.

Case 2: $n=15$ and $X/(g)$ has three singularities of type
$A_{14}$, $A_{4}$ and $A_{2}$.

In any of these cases the discriminant of the lattice $\calR_G$ defined in the previous lemma is a square. So, these cases cannot happen.
\end{proof}

\section{The main theorem}

A Mathieu representation of a finite group $G$ is a 24-dimensional
representation  on a vector space $V$ over a field of characteristic
zero with character
$$\chi(g)=\epsilon({\rm ord}(g)),$$ where
\begin{equation}\label{mumu}
\epsilon(n)=24(n\prod_{p|n}(1+{1\over p}))^{-1}.
\end{equation}
The number
\begin{equation}\label{mu}
\mu(G) = \frac{1}{\# G}\sum_{g\in G}\epsilon(\ord(g))
\end{equation}
is equal to the dimension of the subspace $V^G$ of $V$.
The natural action of a finite group $G$ of symplectic
automorphisms of a complex K3 surface on the singular cohomology
$$H^*(X,{\bbQ})=\oplus_{i=0}^4H^i(X,{\bbQ})\cong {\bbQ}^{24}$$ is a Mathieu representation with
$$\mu(G) = \dim H^*(X,{\bbQ})^G \ge 5.$$ From this Mukai deduces that $G$ is
isomorphic to a subgroup of $M_{23}$ with at least 5 orbits. In
positive characteristic the formula for the number of fixed points
is no longer true and the representation of $G$ on the $l$-adic
cohomology, $l\ne p$,
$$H^*_{\rm et}(X,{\bbQ}_l)=\bigoplus_{i=0}^4H^i_{\rm et}(X,{\bbQ}_l)\cong
{\bbQ}_l^{24}$$ is not Mathieu in general. In this section, using
Theorem \ref{theorem2}, we will show that if $G$ is tame, i.e. the order of $G$
is coprime to $p$, the natural representation of $G$ on
$H^*_{\rm et}(X,{\bbQ}_l)\cong {\bbQ}_l^{24}$ is Mathieu. We will also
show that  $\dim H^*_{\rm et}(X,{\bbQ}_l)^G \ge 5$ under the additional
assumption that either $X$ or a minimal model of $X/G$ is not
supersingular with Artin invariant $\sigma=1$.

First let us recall that the analog of the lattice of transcendental
cycles on a surface $X$ in characteristic $p >0 $ is the group
$T_l\Br(X)$ equal to the projective limit of groups $\Br(X)[l^n]$,
where $\Br(X) = H^2_{\rm et}(X,\bbG_m)$ is the cohomological Brauer group.
Recall that the Kummer sequence in \'etale cohomology \cite{Milne}
gives the exact sequence
\begin{equation}\label{kumm}
0\to \Pic(X)/l^n\Pic(X) \to H_{\rm et}^2(X,\mu_{l^n}) \to \Br(X)[l^n] \to 0.
\end{equation}
Passing to the projective limit we have the exact sequence of $\bbZ_l$-modules
\begin{equation}\label{trans2}
0\to \Pic(X)\otimes \bbZ_l \to  H_{\rm et}^2(X,\bbZ_l) \to
T_l\Br(X) \to 0.
\end{equation}
Tensoring with $\bbQ_l$, we get an exact sequence of
$\bbQ_l$-vector spaces
\begin{equation}\label{trans1}
0\to \Pic(X)\otimes \bbQ_l \to  H_{\rm et}^2(X,\bbQ_l) \to
V_l\Br(X)\to 0.
\end{equation}
It gives the analog of the usual formula for the second Betti number of a surface
$$b_2(X) = \rho(X)+\lambda(X),$$
where $\rho(X)$ is the Picard number of $X$ and $\lambda(X) =
\dim_{\bbQ_l}V_l\Br(X) $ is the Lefschetz number of $X$. Since
$b_2(X)$ and $\rho(X)$ do not depend on $l\ne p$, the Lefschetz
number $\lambda(X)$ does not depend on $l$ either.

\begin{proposition} \label{mainprop} Let $G$ be a finite group of
symplectic automorphisms of a $K3$ surface $X$ defined in
characteristic $p>0$. Assume that $G$ is tame. i.e. the order of $G$ is coprime to $p$.  Then
for any prime $l\ne p$, the natural representation of $G$ on
$H^*_{\rm et}(X,{\bbQ}_l)\cong {\bbQ}_l^{24}$ is Mathieu.
\end{proposition}

\begin{proof} By Theorem \ref{theorem2},
$${\rm ord}(g) \in \{1,\ldots,8\} $$
for all $g\in G$. By Lefschetz fixed point formula,  the character
$\chi(g)$ of the  representation on the $l$-adic cohomology is
equal to the number of fixed points of $g$, which is equal to
$\epsilon({\rm ord}(g))$. This proves the assertion. This also
shows that the action of $G$ on $H^2_{\rm et}(X,{\bbQ}_l)$ is
faithful.
\end{proof}

\begin{lemma}\label{b-1}  Let $G$ be a finite tame group of symplectic
  automorphisms of a K3-surface $X$.
If a nonsingular minimal model $Y$ of $X/G$ is not supersingular, then
\begin{equation}\label{bbb}
\dim H^*_{\rm et}(X,{\bbQ}_l)^G \ge 5.
\end{equation}
\end{lemma}

\begin{proof}
Since $Y$ is not supersingular, $\lambda(Y) > 0$. If $\lambda(Y) =
1$, the Picard number of $Y$ is equal to 21, and hence $Y$ admits
an elliptic fibration. By Artin \cite{Artin3}, the height $h$ of
the formal Brauer group of an elliptic non-supersingular surface
is finite and
$$\lambda(Y) =22-\rho(Y)\ge 2h\ge 2.$$
 Choose
$l$ coprime to the order of $G$.
It is known that $\dim_{\bbQ_l} (V_l\Br(X))^G = \lambda(Y)$
(\cite{Shioda1}, Proposition 5). Taking a $G$-invariant ample
divisor class defining a nonzero element in $(\Pic(X)\otimes
\bbQ_l)^G$, we see that $\dim H_{\rm et}^2(X,\bbQ_l)^G \ge 3$.
Since the characteristic polynomial does not depend on $l\ne p$, this
is true for all $l\ne p$. Together with $H_{\rm et}^0(X,\bbQ_l)$ and
$H_{\rm et}^4(X,\bbQ_l)$ we get \eqref{bbb}.
\end{proof}

It remains to consider the case when $X/G$ is birationally isomorphic to a supersingular K3-surface.
\begin{lemma}\label{crucial} Assume $p\ne 2$. Assume that a K3 surface $X$ admits a symplectic automorphism $g$
of order $2$. Then $X$ admits an elliptic fibration.
\end{lemma}

\begin{proof} As is well-known it suffices to show that $\rho(X) \ge 5$, or, equivalently, $\lambda(X) \le 17$.    Let $U$ be the open set where $g$ acts freely. We know that $X\setminus U$ consists of 8 fixed points of $g$. Let $G = (g), V = U/G.$ We shall use the two spectral sequences employed in the proof of Lemma \ref{samuel}. It is easy to see that they  give the following exact sequence
\begin{equation}\label{spectral}
0\to H^1(G,\Pic(U)) \to \Br(V) \to \Br(U)^{G} \to H^2(G,\Pic(U)).
\end{equation}
Let $t_+$ (resp. $t_-$) be the rank of $g$-invariant (resp. $g$-antiinvariant) part of  $\Pic(U) \cong \Pic(X)$. We have
$$H^1(G,\Pic(U)) = H^1(G,\Pic(X)) = \Ker(1+g^*)/\textup{Im}(1-g^*) \cong (\bbZ/2\bbZ)^{t_-},$$
 $$H^2(G,\Pic(U)) = H^2(G,\Pic(X)) = \Ker(1-g^*)/\textup{Im}(1+g^*) \cong (\bbZ/2\bbZ)^{t_+}.$$
Splitting \eqref{spectral} in two short exact sequences and passing to 2-torsion subgroups we get the following exact sequences of 2-elementary groups
$$0\to (\bbZ/2\bbZ)^{t_-} \to \Br(V)[2] \to A \to (\bbZ/2\bbZ)^{t_-},$$
$$0\to A \to \Br(U)^{G}[2] \to (\bbZ/2\bbZ)^{t} \to 0,$$
where $t\le t_+$. This gives
\begin{equation}\label{ineqx}
\dim_{\bbF_2} \Br(U)^{G}[2] \le \dim_{\bbF_2} A+t\le (\dim_{\bbF_2} \Br(V)[2]-t_-)+t_- +t $$
$$\le \dim_{\bbF_2} \Br(V)[2]+t_+.
\end{equation}
Let $Y$ be a minimal resolution of singularities of $X/G$ and $\calE$ be the exceptional divisor. According to \cite{DeMF}, the exact sequence of local cohomology for the pair $(Y,\calE)$ and the sheaf $\bbG_m$ defines an exact sequence (modulo $p$-groups)
$$0\to \Br(Y) \to \Br(Y\setminus \calE) \to H^1(\calE,\bbQ/\bbZ).$$
Since $\calE$ is the disjoint union of 8 smooth rational curves, we obtain
$$\Br(Y) \cong \Br(Y\setminus \calE) \cong \Br(V).$$
Similarly, we obtain
$$\Br(U) \cong \Br(X).$$
It follows from \eqref{trans2} that, up to a finite group,
$$\Br(Y)[2] = (\bbZ/2\bbZ)^{\lambda(Y)}.$$
Applying \eqref{ineqx}, we obtain
\begin{equation}\label{ineqxx}
\dim_{\bbF_2} \Br(X)^{G}[2]  \le \lambda(Y)+t_+.\end{equation}
If $\rho(X) \ge 5$ we are done. Otherwise $\rho(X) \le 4$, hence $t_+ \le 4$. Since $Y$ contains 8 disjoint smooth rational curves and also the pre-image of a  class of an ample divisor on $X/G$, we have $\rho(Y) \ge 9$, and therefore $\lambda(Y) \le 22-9 = 13$. Now \eqref{ineqxx}  implies
\begin{equation}\label{ineqxxx}
\dim_{\bbF_2} \Br(X)^{G}[2] \le 17.
\end{equation}
The exact sequence of sheaves in \'etale topology
$$0\to \mu_{2^{n}} \overset{[2]}{\longrightarrow} \mu_{2^{n+1}}\to \mu_2 \to 0$$
gives, after passing to cohomology and taking the projective limits, the exact sequence
$$H^2(X,\bbZ_2) \overset{[2]}{\longrightarrow} H^2(X,\bbZ_2) \to H^2(X,\mu_2) \to 0.$$
Since an automorphism  of order $\le 2$ of a free $\bbZ_l$-module acts trivially modulo 2, we obtain that
$$ H^2(X,\mu_2)^G = H^2(X,\mu_2).$$
Applying the Kummer exact sequence \eqref{kumm}, we see that
$$\Br(X)^{G}[2] = \Br(X)[2].$$
It remains to apply  \eqref{ineqxxx}.
\end{proof}

\begin{lemma}\label{sup} Let $G$ be a finite tame group of symplectic
automorphisms of a K3-surface $X$ of order  $\ne 7, 21$. Assume
that a minimal resolution $Y$ of $X/G$ is a supersingular K3
surface. Then $X$ is supersingular.
\end{lemma}

\begin{proof} Recall from \cite{Artin3} that the formal Brauer group
  $\hat{\Br}(S)$ of a supersingular K3-surface $S$  is isomorphic to
  the formal additive group $\hat{\bbG}_a$. It is conjectured that the
  converse is true, and it has been verified  if $S$ is an elliptic
  surface (loc.cit. Theorem (1.7)). Since Brauer group is a birational invariant, the projection
$\pi: X\to X/G$ defines a natural homomorphism of the formal Brauer
groups
$\pi^*:\hat{\Br}(Y)\to \hat{\Br}(X).$ The corresponding map of the tangent spaces is
$\pi^*:H^2(Y,\calO_Y)\to H^2(X,\calO_X)$. Since the order of $G$ is
prime to the characteristic, the trace map shows that this
homomorphism is nonzero. Since there are no non-trivial maps between a
formal group of finite height and $\hat{\bbG}_a$ we obtain that
$\hat{\Br}(X)\cong \hat{\bbG}_a.$ If $G$ contains an element of order
2, we are done by Lemma \ref{crucial}. Assume that $G$ has no elements
of order 2. By
Proposition \ref{mainprop}, the representation of $G$ in
$H^*(X,\bbQ_l)$ is a Mathieu representation.
Following Mukai's arguments from \cite{Mukai}, we obtain that the
order of $G$ must divide $3^2.5.7$.

Suppose that 3 distinct prime numbers divide $\# G$. It is known that
no simple non-abelian group of order dividing $3^2.5.7$ exists. Thus
$G$ is solvable and hence contains a subgroup of order $35$
(\cite{Hall}, Theorem 9.3.1). It follows easily from Sylow's Theorem
that such a group is cyclic and hence is not realized as a group of
symplectic automorphisms of $X$.

Suppose  $\#G= 3^a.5$ or $3^a.7$ with $a\ne 0$.  Again, by Sylow's theorem we obtain that a Sylow 5-subgroup (or 7-subgroup)  is normal. Since no element of order 3 can commute with an element of order 5 or 7, we obtain that $G$ is a non-abelian group of order 21. This case is excluded by the assumption.

The remaining possible cases are $\#G = 3, 5, 9$. This gives that
$X/G$ has either 6 singular points of type $A_2$, or $4$ singular
points of type $A_4$, or 8 singular points of type $A_2$.
Since $\rho(Y) = 22$, this immediately implies that $\dim \Pic(X/G)\otimes \bbQ = 22-12 =10$ or $22-16 = 6$. Hence $\rank \Pic(X)^G \ge 6$. Thus $X$ is an elliptic surface, and, by Artin's result cited above, we obtain that $X$ is supersingular.
\end{proof}

\begin{proposition}\label{mainmain} Let $G$ be a finite tame group of
  symplectic automorphisms of a K3 surface $X$.
  Assume  that either $X$ or a minimal model of $X/G$ is
  not a supersingular K3 surface with Artin invariant $\sigma = 1$.
Then
$$\dim H_{\et}^*(X,\bbQ_l)^G \ge 5.$$
\end{proposition}

\begin{proof} By Lemma \ref{b-1} we
  may assume that a minimal nonsingular model $Y$ of $X/G$ is
  supersingular. A symplectic group of order 7 or 21 is uniquely
  determined and is in Mukai's list \cite{Xiao} and satisfies
the assertion of the proposition. By Lemma \ref{sup} we obtain that $X$ is supersingular.

Assume that the Artin invariant $\sigma$ of $X$ is greater than 1.

Let us  consider the representation of $G$
on the crystalline cohomology $H_{\crys}^*(X/W)$. We refer to
\cite{Illusie} for the main properties of crystalline cohomology
and to \cite{RS}, \cite{Ogus} for particular properties of
crystalline cohomology of K3 surfaces. The cohomology
$H_{\crys}^*(X/W)$, where $X$ is a K3-surface, is a free module of
rank 24 over the  ring of Witt vectors $W = W(k)$. The vector
space $H_{\crys}^*(X/W)_K$, where $K$ is the field of fractions of
$W$, is of dimension 24. The ring $W$ is a complete noetherian
local ring of characteristic 0 with maximal ideal $(p)$ and the
residue field isomorphic to $k$. The quotient module
$H_{\crys}^*(X/W)/pH_{\crys}^*(X/W)$ is a $k$-vector space of
dimension 24 isomorphic to the algebraic De Rham cohomology
$H_{\DR}^*(X)$. Let
$$H = H_{\DR}^2(X).$$ It is known  that the Hodge
spectral sequence
$$E_1^{p,q} = H^q(X,\Omega_X^p)\Rightarrow H_{\DR}^n(X)$$
degenerates and we have the following canonical exact sequences:
\begin{eqnarray}\label{}
\label{ }
0 \to F^1H\to H \to  H^2(X,\calO_X) \to 0,\\
0\to H^0(X,\Omega_X^2)\to F^1H \to H^1(X,\Omega_X^1) \to 0.
\end{eqnarray}
Here
$$F^iH = \sum_{p\ge i} H^{q}(X,\Omega_X^p)\cap H$$ is the
Hodge filtration of the De Rham cohomology. Obviously, the
subspace $F^1H$ is $G$-invariant. Since the order of  $G$ is prime
to the characteristic, the representation of $G$ is semi-simple
and  hence the $G$-module $H$ is isomorphic to the direct sum of
$G$-modules
\begin{equation}\label{sum}
H \cong H^0(X,\Omega_X^2)\oplus H^2(X,\calO_X)\oplus H^1(X,\Omega_X^1).
\end{equation}
By definition,
$$H^0(X,\Omega_X^2) \subset H^G.$$ By Serre's duality
$$H^2(X,\calO_X)\subset H^G.$$
This shows that $\dim H^G \ge 2$.

Let
$$V = H^2_{\crys}(X/W).$$ The multiplication by $p$ map $[p]$
defines the exact sequence
$$0\to V\overset{[p]}{\to} V \to H\to 0.$$
Taking $G$-invariants we obtain the exact sequence
\begin{equation}\label{exxact}
0\to V^G\overset{[p]}{\to} V^G \to H^G \to H^1(G,V).
\end{equation}
Since the ring $W$ has characteristic 0, the multiplication by
$|G|$ defines  an injective map $V\to V$. On the other hand, it
induces the zero map on the cohomology $ H^1(G,V)$. This implies
that
$$H^1(G,V) = 0.$$ This shows that $V^G$ is a free submodule of
$V$ of rank equal to $\dim H^G $.

 It is known  that the Chern class map
$$c_1:S_X\to H_{\crys}^2(X/W)$$ is injective and its composition with the
reduction mod $p$ map
$$H_{\crys}^2(X/W)\to H^2_{DR}(X)$$ defines an injective map
\begin{equation}\label{chern}c: S_X/pS_X \to H^2_{DR}(X)
\end{equation}
with image contained in $F^1H^2_{DR}(X)$ (see \cite{Ogus}).

If $X$ is supersingular with Artin invariant $\sigma > 1$, the
composition of $c$ with the projection $F^1H^2_{DR}(X) \to
H^1(X,\Omega_X^1)$ is injective. This result is implicitly contained
in \cite{Ogus} (use Remark 2.7 together with the fact that a
supersingular surface with $\sigma > 1$ admits a non-trivial
deformation to a supersingular surface with Artin invariant equal to
1). Let $L$ be a $G$-invariant ample line bundle on $X$. We may assume
that its isomorphism class defines a non-zero element in
$S_X/pS_X$. Thus its image in $H^1(X,\Omega_X^1)$ is a nonzero
$G$-invariant element and we get 3 linearly independent elements in
$H^G$, each from one of the three direct summands of $H$. Applying
\eqref{exxact}, we find 3 linearly independent elements in
$H_{\crys}^2(X/W)^G$.

Since $H_{\crys}^0(X/W)$,  $H_{\crys}^4(X/W)$ are trivial
$G$-modules, we obtain
\begin{equation}
\label{bound5} \dim H_{\crys}^*(X/W)_K^G \ge 5.
\end{equation}
It remains to use  the fact that  the characteristic polynomials
of $g\in G$ on  $ H_{\crys}^*(X/W)_K$ and on $ H_{\et}^*(X,\bbQ_l)$,
$l\ne p$, have integer coefficients and coincide with each other
(\cite{Illusie}, 3.7.3).

Thus if the assertion is not true, $X$ must be supersingular of Artin invariant $\sigma = 1$.

 Assume that $Y$ is supersingular with Artin invariant $\sigma > 1$.  Let $X'$ be the open subset of $X$ where $G$ acts freely and let $Y' = X'/G$. The standard Hochshild-Serre spectral sequence  implies that the pull-back under the projection $f:X'\to Y'$ defines an isomorphism
$$\Pic(Y')/\Hom(G,k^*) \cong \Pic(X')^G \cong \Pic(X)^G.$$
Let $R$ be the sublattice of $S_Y$ spanned by the  irreducible components of exceptional curves of the resolution $\pi:Y\to X/G$. It is isomorphic to the orthogonal sum of root lattices of discriminants prime to $p$. The restriction map $\Pic(Y)\to \Pic(Y')$ is surjective and its kernel is $R$. The torsion group of $\Pic(Y')$ is isomorphic to $\Hom(G,k^*)$. Let $R'$ be the saturation of $R$ in $\Pic(Y)$. We have
\begin{equation}\label{lat}
N: = \Pic(Y)/R' \cong \Pic(Y')/\Hom(G,k^*) \cong \Pic(X)^G.
\end{equation}
Since $\#G$ is coprime to $p$, the discriminant of the sublattice $R'$
 is coprime to $p$.  The discriminant group $D_Y$ of $\Pic(Y)$ is an
 elementary $p$-group of rank $2\sigma \ge 4$ and is a subquotient of
 the discriminant group of $R'\oplus R^\perp$. This implies that
 $\rank N = \rank R^\perp > 2$,  it follows from \eqref{lat} that
 $\rank \Pic(X)^G > 2$, and  by the
 above arguments we will find 5 linearly independent elements in $H_{crys}^*(X/W)^G$.
\end{proof}

\begin{remark} If $p$ divides $|G|$, the exact sequences (4.11) and
(4.12) may not split as $G$-modules. In fact, there are examples
where $$\dim H_{\crys}^*(X/W)_K^G = 4,$$ so that \eqref{bound5}
does not hold.
\end{remark}

\begin{theorem}\label{thm4.6} Let $G$ be a finite group of symplectic
automorphisms of a K3 surface $X$. Assume that $G$ is tame and
that either $X$ or a minimal model of $X/G$ is not a supersingular
K3 surface with Artin invariant $\sigma = 1$. Then $G$ is a
subgroup of the Mathieu group $M_{23}$ which has $\ge 5$ orbits in
its natural permutation action on the set of $24$ elements. All
such groups are subgroups of the $11$ groups listed in
\cite{Mukai}.
\end{theorem}

\begin{proof} Let us  consider the linear representation $\rho$ of
$G$ on  $H_{\et}^*(X,\bbQ_l), l\ne p$.  Applying Proposition
\ref{mainprop} and
\ref{mainmain}, we find that $\rho$ is a Mathieu representation
over the field $\bbQ_l$ with $\dim H_{\et}^*(X,\bbQ_l)^G \ge 5$.
Replacing $\bbQ$ with $\bbQ_l$ we repeat
the arguments of Mukai. He uses at several places the fact that
the representation is over $\bbQ$.  The only essential place where
he uses that the representation is over $\bbQ$ is Proposition
(3.21), where $G$ is assumed to be a 2-group containing a maximal
normal abelian subgroup $A$ and the case of $A =(\bbZ/4)^2$ with
$\#(G/A)\ge 2^4$ is excluded by using that a certain 2-dimensional representation of
the quaternion group $Q_8$ cannot be defined over $\bbQ$. We use that $G$
also admits a Mathieu representation on 2-adic cohomology, and it
is easy to see that the representation of  $Q_8$ cannot be defined
over $\bbQ_2$. To show that the 2-Sylow subgroup of $G$ can be
embedded in $M_{23}$ he uses the fact that the stabilizer of any
point on $X$ is isomorphic to a finite subgroup of $\SL(2,\bbC)$,
and the classification of such subgroups allows him  to exclude
some groups of order $2^n$. By Lemma \ref{klein}, we have the same classification, so we can do the same.
\end{proof}

Applying  Theorem \ref{main}, we obtain the following.

\begin{corollary} Assume that $p > 11$ and  either $X$ or a minimal model of $X/G$ is not a supersingular K3
surface with Artin invariant $\sigma = 1$. Then $G$ is a subgroup of $M_{23}$ with $\ge 5$ orbits and hence belongs to Mukai's list.
\end{corollary}

\section{The exceptional case}

Here we investigate the case when the order of $G$ is prime to $p$ and
\begin{equation}\label{mu4}
\dim H_{\crys}^*(X/W)_K^G = 4.
\end{equation}
By theorem \ref{thm4.6} this may happen only if   both $X$ and a minimal
nonsingular model $Y$ of $X/G$ is a supersingular K3 surface with
Artin invariant $\sigma = 1$.
We  refer to this as the exceptional case and the group $G$ will be
called an
\emph{exceptional group}.

It is known that a supersingular surface with Artin invariant $\sigma = 1$ is unique up to isomorphism. More precisely, we have the following (see \cite{Ogus},Corollary 7.14).

\begin{proposition}\label{ogus}
Let $X$ be a supersingular surface with Artin invariant $\sigma = 1$. Assume that $p \ne  2$. Then $X$ is birationally isomorphic to  the Kummer surface of the abelian surface $E\times E$, where $E$ is a supersingular elliptic curve.
\end{proposition}

If $p = 2$, the surface is explicitly described in \cite{DKo}.  Note that the Kummer surface does not depend on $E$. If $p\equiv 3 \mod 4$ (resp. $p\equiv 2 \mod 3$) we can take for $E$ an elliptic curve with Weierstrass equation $y^2 = x^3-x$ (resp. $y^2 = x^3+1$).

\bigskip
 It follows from the proof of
Proposition \ref{mainmain} and Lemma \ref{square} that an exceptional group satisfies the following properties:
\begin{itemize}
\item[(EG1)] $G$ admits a Mathieu representation $V_l$ over any $\bbQ_l, l\ne p$;
\item[(EG2)] $\mu(G) = \dim V_l^G = 4$;
\item[(EG3)] the root lattice $\calR_G$ spanned by irreducible
  components of the exceptional locus of the resolution $Y\to X/G$ is
  of rank 20 (this is equivalent to (EG2));
\item[(EG4)] $d_{\calR_G}$ is coprime to $p$ and is not a square.

\end{itemize}

\bigskip
We use the following notations of groups from \cite{Mukai} and the ATLAS \cite{ATLAS}:

\bigskip
 $C_n$ the cyclic group of order $n$, sometimes denoted by $n$,

 $D_{2n}$ the dihedral group of order $2n$,

 $Q_{4n}$ the binary dihedral group of order $4n$,

 $T_{24}$ the binary tetrahedral group,

 $O_{48}$ the binary octahedral group,

$\frakS_n$ the symmetric group of degree $n$,

 $\frakA_n$ the alternating group of degree $n$,

$\frakS_{n_1,\ldots,n_k}$ a subgroup of  of $\frakS_{n_1+\ldots+n_k}$
which preserves the decomposition of a set of
  $n_1+\ldots+n_k$ elements as a disjoint union of subsets of cardinalities $n_1,\ldots,n_k$.

 $\frakA_{n_1,\ldots,n_k}$ = $\frakS_{n_1,\ldots,n_k}\cap \frakA_{n_1+\ldots+n_k}$;

 $M_{k}$ the Mathieu group of degree $k$.

$L_n(q) = \textup{PSL}_n(\bbF_q)$.

$A_{\text{\Tiny \textbullet}} B$ a group which has a normal subgroup isomorphic to $A$ with quotient isomorphic to $B$;

$A^{\text{\Tiny \textbullet}} B$ as above but the extension does not split;

$A:B$ a semidirect product with normal subgroup $A$;

$A\circ B$ the central product of two groups.

We will also use the notations for 2-groups from \cite{HS} and \cite{Mukai}.

\bigskip
The goal of this section is to prove the following.

\begin{theorem}\label{except} An exceptional group $G$  is
  isomorphic to one of the following groups (the corresponding
root lattice  $\calR_G$ is given  in the parenthesis):
\begin{enumerate}
\item[(I)] Non-solvable groups:
\begin{enumerate}
\item[(i)]  $\frakS_6$ of order $2^4.3^2.5$ $(A_4+2A_3+2A_5)$;
\item[(ii)]  $M_{10}$ of order $2^4.3^2.5$ $(A_4+A_3+A_2+A_7+D_4)$;
\item[(iii)] $2^4:\frakA_6$ of order $2^7.3^2.5$ $(E_7+A_2+A_3+2A_4)$;
\item[(iv)] $M_{21}$ of order $2^6.3^2.5.7 \ (A_2+2A_4+A_6+D_4)$;
\item[(v)] $\frakA_7$ of order $2^3.3^2.5.7\ (A_2+A_3+A_4+A_6+D_5)$;
\item[(vi)] $M_{20}' = 2^4:\frakA_5(\ncong M_{20})$ of order $2^6.3.5 \ (A_1+2A_4+D_5+E_6)$;
\item[(vii)] $M_{20}:2\cong 2^4:\frakS_5$ of order $2^7.3.5 \ (A_2+A_3+A_4+A_5+D_6)$;
\item[(viii)] $2^3:L_2(7)$ of order $2^6.3.7 \ (A_2+2A_3+A_6+E_6)$.
\end{enumerate}
\end{enumerate}
\begin{enumerate}
\item[(II)] Solvable groups:
\begin{enumerate}
 \item[(ix)] $3^2:C_8$ of  order $2^3.3^2$  $( 2A_7+A_3+A_2+A_1)$;
\item[(x)] $3^2:SD_{16}$ of order $2^4.3^2$
   $( A_7+A_5+D_4+A_3+A_1)$;
\item[(xi)]  $2^3: 7$ of order $2^3.7$
  $ (3A_6+2A_1)$;
\item[(xii)] $2^4:(5:4)$ of order $2^6.5$
   $(A_7+A_4+3A_3)$;
\item[(xiii)]$2^2\,_{\text{\Tiny \textbullet}} (\frakA_4\times \frakA_4)=\Gamma_{13}a_1:3^2$
of order $2^6.3^2$  $(E_6+2A_2+2A_5)$;
\item[(xiv)] $2^2\,_{\text{\Tiny \textbullet}} \frakA_{4,4}=\Gamma_{13}a_1:\frakA_{3,3}=2^4:\frakA_{3,4}$ of order
$2^7.3^2$ $(E_7+2D_5+A_2+A_1)$;
\item[(xv)] $2^4:(3\times D_6)$ of order $2^5.3^2$
   $(3A_5+A_3+A_2)$;
\item[(xvi)] $2^4:(3^2: 4)$ of order $2^6.3^2$
   $( A_7+D_5+2A_3+A_2 )$;
\item[(xvii)] $2^4:\frakS_{3,3}=\frakS_{4,4}=2^4:\frakA_{2,3,3}$ of order $2^6.3^2$
   $(D_5+D_4+2A_5+A_1)$;
\item[(xviii)]  $O_{48}$ of order $2^4.3$ $(E_7+D_6+D_5+A_2$ or  $2E_7+D_4+A_2)$;
\item[(xix)] $T_{24}\times 2$ of order $2^4.3$ $( E_6+D_4+2A_5)$;
\item[(xx)] $O_{48}:2$ of order $2^5.3$ $(E_7+D_6+A_5+2A_1)$;
\item[(xxi)] $(Q_8\circ Q_8)^{\text{\Tiny \textbullet}}\frakS_3=
\Gamma_5a_1\, ^{\text{\Tiny \textbullet}}\frakS_3$  of order $2^6.3$ $(A_1+A_2+A_3+2E_7)$;
\item[(xxii)] $(2^4:2)^{\text{\Tiny \textbullet}}\frakS_3=2^4:Q_{12}$  of order $2^6.3$ $(A_1+A_2+A_7+2D_5)$;
\item[(xxiii)] $\Gamma_{13}a_1:3=2^4:\frakA_4$  of order $2^6.3$ $(A_1+3A_5+D_4$ or  $2A_1+3E_6)$;
\item[(xxiv)] $\Gamma_{13}a_1:\frakS_3=2^4:\frakS_4$  of order $2^7.3$
 $(A_1+A_3+A_5+D_5+D_6)$;
\item[(xxiv$'$)] $\Gamma_{13}a_1:\frakS_3=2^4:\frakS_4$  of order $2^7.3$
 $(A_1+2A_3+E_6+E_7)$;
\item[(xxv)] $\Gamma_{25}a_1:3=2^4:\frakA_4$  of order $2^6.3$ $(A_1+A_3+2A_5+E_6)$;
\item[(xxvi)] $\Gamma_{25}a_1:\frakS_3=2^4:\frakS_4=\Gamma_5a_1\, ^{\text{\Tiny \textbullet}}D_{12}$  of order $2^7.3$
  $(A_1+A_3+A_5+D_4+E_7)$;
\item[(xxvii)] $(2^3:7):3$ of order $2^3.3.7$
$(A_6+2A_5+2A_2)$.
\end{enumerate}
\end{enumerate}
\end{theorem}

\begin{remark}
(1) All these groups are subgroups of the Mathieu group $M_{23}$
with number of orbits equal to 4 (see Proposition \ref{subM23}).
Note that  not all subgroups of $M_{23}$ with $4$ orbits belong to our list.  For example, a group containing elements of order $>8$ must be wild  by Theorem
\ref{theorem2}, hence is not contained in our list.  Examples of such groups are $L_2(11), 2\times L_2(7), \frakA_{3,5}$.  We thank D. Allcock and S. Kond\=o for confirming
this. Also, as we will see in the proof of Lemma 5.10, Case 3,  there is a degree 2 extension of $M_{20}'$ isomorphic to $2^4:S_5$ which cannot act symplectically.  According to Allcock, this group can be realized as a subgroup of $M_{23}$ preserving the partition $(5,2,1,16)$ with $(5,2,1)$ forming an octad of the Steiner system. Our double extension $M_{20}:2 \cong 2^4:S_5$ preserves the partition $(1,1,2,20)$. 

(2) Most exceptional groups seem to realize in infinitely many
different characteristics. A full account on the realization
problem will be given in other publication.

(3) Among the 27 groups the following 10 groups are maximal:
(i)-(v), (vii), (viii), (x), (xiv), (xvii).
\end{remark}

We will prove the theorem by analyzing and extending Mukai's arguments from
\cite{Mukai}.
We can use the arguments only based on property (EG1) of $G$ and do
not use the assumption that $\mu(G) \ge 5$.

Recall  that
 $$\mu(G) = \frac{1}{\#G}\sum_{g\in G}\epsilon(\ord(g)),$$ where
 $\epsilon(n)$ is given in  \eqref{mumu}. Let $g$ be an element of order $n$ prime to $p$ acting symplectically on a K3 surface. Using Theorem \ref{theorem2} and computation of
 $\#X^g$ from \cite{Ni1}, one checks that
\begin{equation}\label{fix}
\epsilon(n) = \#X^g.
\end{equation}
The possible values of $\epsilon(n)$ are given in Theorem \ref{theorem2}.

\begin{lemma}\label{xiao}
$$\sum_{Gx\in X/G} \frac{1}{\#G_x} = \frac{24}{\#G}+k-\mu(G),$$
where $k$ is the number of singularities on $X/G$.
\end{lemma}

\begin{proof}  Let
$$S = \{(x,g)\in X\times G\setminus \{1\}: g(x) = x.\}$$
By projecting to $X$ we get
$$\#S = \#G\sum_{Gx\in \textup{Sing}(X/G)}(1-\frac{1}{\#G_x}).$$
By projecting to $G$, and using \eqref{fix}, we get
$$\#S = \sum_{g\in G\setminus \{1\}}\#X^g =  \sum_{g\in G\setminus \{1\}}\epsilon(\ord(g))  =\#G\mu(G)-24.$$
This gives
 $$\sum_{Gx\in \textup{Sing}(X/G)} \frac{1}{\#G_x} = \frac{24}{\#G} + k-\mu(G).$$
\end{proof}

\begin{remark}
We also have the following formula from \cite{Xiao}
$$\sum_{Gx\in X/G}\frac{1}{\#G_x} = \frac{24}{\#G}+k+\rank R_G-24.$$
Comparing this with the previous formula, we get
\begin{equation}\label{xiao2}
\mu(G) = 24-\rank R_G.
\end{equation}
\end{remark}

The classification of finite subgroups of $\SL(2,k)$ which admit a
Mathieu representation is given in \cite{Mukai}, Proposition
(3.12). The groups $Q_{4n}$, $n\ge 5$, and the binary icosahedral
group $I_{120}$ are not realized since they contain an element of
order $>8$. The following table gives the information about possible
stabilizer groups $G_x$, their orders $o_x$, the number $c_x$ of
irreducible components in a minimal resolution, types of singular
points and the structure of the discriminant group $D_x$ of the corresponding root lattice.

\begin{table}[htdp]
\begin{center}
\begin{tabular}{|l||r|r|r|r|r|r|r|r|r|r|r|r|}
\hline
$G_x$&$C_2$&$C_3$&$C_4$&$C_5$&$C_6$&$C_7$&$C_8$&$Q_8$&$Q_{12}$&$Q_{16}$&$T_{24}$&$O_{48}$\\ \hline
$o_x$&2&3&4&5&6&7&8&8&12&16&24&48\\ \hline
$c_x$&1&2&3&4&5&6&7&4&5&6&6&7\\ \hline
Type&$A_1$&$A_2$&$A_3$&$A_4$&$A_5$&$A_6$&$A_7$&$D_4$&$D_5$&$D_6$&$E_6$&$E_7$\\ \hline
$D_x$&$C_2$&$C_3$&$C_4$&$C_5$&$C_6$&$C_7$&$C_8$&$C_2^2$&$C_4$&$C_2^2$&$C_3$&$C_2$\\ \hline
\end{tabular}
\end{center}
\caption{}
\end{table}

Note that $\mu(O_{48}) = 4$ so this case may occur only in the exceptional
case.

For an  exceptional group $G$, using the above table and the formula
from Lemma \ref{xiao} it is easy to show that
\begin{itemize}
\item[(EG5)] the number $k$ of singularities on $X/G$ is 4 or 5.
\end{itemize}

Let $G$ be an exceptional group acting on $X$. It defines a set of $k$ numbers
$o_1,\ldots, o_k$ and $c_1,\ldots,c_k$ from  Table 1. We have
\begin{itemize}
\item[(i)] $c_1+\ldots+c_k = 20$;
\item[(ii)] $\frac{1}{o_1}+\ldots+\frac{1}{o_k} = k-4+\frac{24}{N}$,
 where $N$ is a positive integer (equal to $\#G$);
\item[(iii)] $o_i| N$ for all $i = 1,\ldots,k$;
\item[(iv)] $d_1\dots d_k$ is not a square, where $d_i$ is the order
  of the discriminant group $D_i$;
\item [(v)] $k = 4$ or $5$.
\end{itemize}

\bigskip
We are grateful to Daniel Allcock who run for  us a computer program
which enumerates all collections of numbers $o_1,\ldots,o_k$ and the
corresponding numbers $c_i,d_i,N$ satisfying properties (i)-(v). This
gives us all possible orders $N$ of possible exceptional groups $G$ as
well as all possible root lattices describing a minimal resolution of
 singularities of $X/G$. We refer to this as the {\it List}.  It is reproduced in Table 2.
\begin{table}[h]
\begin{center}
\begin{tabular}{|l||r|r}
\hline
Order&Root lattices\\ \hline
$2^6.3^2.5.7$&$A_2A_4A_4A_6D_4$\\ \hline
$2^3.3^2.5.7$&$A_2A_3A_4A_6D_5$\\ \hline
$2^4.5.7$&$A_3A_3A_4A_4A_6$\\ \hline
$3^2.5.7$&$A_2A_4A_4A_4A_6$\\ \hline
$2^7.3^2.5$&$A_2A_3A_4A_4E_7$\\ \hline
$2^5.3^2.5$&$A_3A_4A_4A_4A_5$\\  \hline
$2^4.3^2.5$&$A_1A_3A_4E_6E_6,\ A_1A_4A_5D_4E_6,\ A_1A_4A_5D_5D_5$\\
&$A_2A_3A_4A_7D_4,\ A_3A_3A_4A_5A_5$\\ \hline
$2^3.3^2.5$&$A_2A_4A_4A_5A_5$\\ \hline
$2^6.3.7$&$A_2A_3A_3A_6E_6,\ A_2A_3A_5A_6D_4$\\ \hline
$2^3.3.7$&$A_1A_3A_5A_6D_5,\ A_2A_2A_5A_5A_6$\\ \hline
$2^7.3.5$&$A_2A_3A_4A_5D_6$\\ \hline
$2^6.3.5$&$A_1A_4A_4D_5E_6$\\ \hline
$2^3.3.5$&$A_1A_2A_4A_7E_6,\ A_1A_4A_5A_5A_5$\\ \hline
$2^7.3^2$&$A_1A_2D_5D_5E_7,\ A_1A_2D_5D_6E_6$\\ \hline
$2^6.3^2$&$A_1A_2D_4D_6E_7, \ A_1A_2D_5D_6D_6,\ A_1A_5A_5D_4D_5$,\\
&$A_2A_2A_2E_7E_7,\ A_2A_2A_5A_5E_6, \ A_2A_3A_3A_7D_5$ \\ \hline
$2^5.3^2$&$A_2A_2A_5A_7D_4,\ A_2A_3A_5A_5A_5,\ A_2A_2A_3A_7E_6$\\ \hline
$2^4.3^2$&$A_1A_1D_5D_6E_7,\ A_1A_1D_6D_6E_6,\ A_1A_3A_5A_7D_4,\ A_1A_3A_3A_7E_6$\\ \hline
$2^3.3^2$&$D_4D_5D_5E_6, \ A_1A_2A_3A_7A_7$\ \\ \hline
$2^3.7$&$A_1A_1A_6A_6A_6$\\ \hline
$2^6.5$&$A_1A_3A_4D_6D_6,\ A_1A_4A_7D_4D_4,\ A_3A_3A_3A_4A_7$\\ \hline
$2^7.3$&$A_1A_3A_3E_6E_7,\ A_1A_3A_5D_4E_7,\ A_1A_3A_5D_5D_6$\\ \hline
$2^6.3$&$A_1A_1E_6E_6E_6,\ A_1A_2A_3E_7E_7,\ A_1A_2A_7D_5D_5$,\\
&$ A_1A_3A_5A_5E_6,\ A_1A_5A_5A_5D_4$\\ \hline
$2^5.3$&$A_1A_1A_5D_6E_7, \ A_1A_1A_7D_5E_6,\ A_2A_2A_2A_7A_7$\\ \hline
$2^4.3$&$A_2D_4E_7E_7, \ A_2D_5D_6E_7, \ A_3A_5E_6E_6, \ A_5A_5D_4E_6$\\ \hline

\end{tabular}
\end{center}
\caption{}
\end{table}

\begin{lemma}\label{generators} Let $l\ne p$ be a prime number. Assume that the minimal number of generators of the $l$-torsion part of the discriminant group  of the lattice $\calR_G$ is greater than 2. Then $l$-part of the abelianized group $G/[G,G]$ is non-trivial.
\end{lemma}

\begin{proof} Let $M$ be the saturation of the lattice $\calR_G$ in
  $\Pic(Y)$. It follows from  the proof of  Proposition \ref{mainmain}
  that $G/[G,G] \cong M/\calR_{G}$. The lattice $M\oplus M^\perp$ is a
  sublattice of finite index of the lattice $S_Y$ with discriminant
  group $(\bbZ/p\bbZ)^2$. Tensoring with the localized ring $\bbZ_{(l)}$
  at the prime ideal $(l)\subset \bbZ$,  we obtain
 $$M_l\oplus (M^\perp)_l \subset (S_Y)_l.$$ Since $(S_Y)_l$ is
 unimodular over $\bbZ_{(l)}$,
 the discriminant groups of $M_l$ and  $(M^\perp)_l$ are isomorphic to
 each other. Since $\rank_{\bbZ_{(l)}}~ (M^\perp)_l =2$, we obtain
that the discriminant group of $M_l$ is generated by
  $\le 2$ elements. If $l$ does not divide
$\# G/[G,G]$, then $M_l \cong (\calR_G)_l$. The discriminant group  of
$(\calR_G)_l$ is the $l$-torsion part of the discriminant group of
$\calR_G$, and the assertion follows.
\end{proof}

\begin{lemma}\label{35} If $G$ is of order divisible by $35$, then $G$
  is isomorphic to the simple group $M_{21} \cong L_3(4)$ of order
  $2^6.3^2.5.7$ or the alternating group $\frakA_7$ of order $2^3.3^2.5.7$.
\end{lemma}

\begin{proof}
Let $[G,G]$ be the commutator subgroup of $G$.

\medskip\noindent
Case 1. $G = [G,G]$.

\medskip
By Lemma \ref{generators}, for any prime $l\ne
p$ the $l$-part of the discriminant group of $\calR_G$ is generated by $\le 2$
elements.

The List
shows that there are 4 possible orders  $N$ divisible by 35. We will
eliminate the orders  $N =  2^4.5.7$ and $3^2.5.7$.  In the first
case, considering a 7-Sylow subgroup, we see that no elements of order 5 or 2 can
normalize it, because $G$ does not contain elements of order $35$, or
14 and does not contain $D_{14}$, which does not admit a Mathieu representation.
By Sylow's theorem, $\#Syl_7(G) = 2^4.5 \equiv 1 \mod 7$, a
contradiction.  In the second case, no elements of order 7 or 3 can
normalize a 5-Sylow subgroup, and hence $\#Syl_5(G) =3^2.7 \equiv
1 \mod 5$, again a contradiction.

Assume $\#G = 2^3.3^2.5.7$. This is the order of $\frakA_7$. Let us
show that $G$ is simple. Assume $G$ is not simple and let $H$ be a
normal subgroup such that $G/H$ is simple (non-abelian because
$G= [G,G]$). It follows from \cite{Mukai}, Proposition (3.3) that
$\mu(G/H) > \mu(G) = 4. $ The group $G/H$ acts symplectically on a
minimal resolution of $X/H$ and belongs to Mukai's list. It follows
from Proposition (4.4) of loc.cit. that $G/H\cong \frakA_5, \frakA_6$,
or $L_2(7)$. In the first case $\#H = 42$. There is no such group in
Mukai's list as well as  in our List. In the second case $\#H = 7$. It is
known that $\frakA_6$ does not admit such a nontrivial extension
(necessarily central). In the last case $\#H = 15$ and again we use that
$G$ does not contain an element of order $15$. Thus $G$ is
simple. It is known that there is only one simple group of order equal to the order of $\frakA_7$. The List gives only
one possible root lattice  $\calR_G=A_2+A_3+A_4+A_6+D_5$.

Assume $\#G = 2^6.3^2.5.7$. This is the order of the Mathieu group
$M_{21}$. As in the previous case we show that $G$ is simple by
analyzing the kernel $H$ of a homomorphism onto a simple quotient of
$G$. As before,  $G/H\cong \frakA_5, \frakA_6$, or $L_2(7)$. In the
first case $\#H = 2^4.3.7$. One checks that there are no  such groups
in Mukai's list (Theorem (5.5) and Proposition (4.4)) and in our
List. In the second case, $\#H = 2^3.7$  and a group of this order
with $\mu(G) = 4$ is possible. We will show later that $H$ must be
isomorphic to $C_2^3:C_7$. The Sylow subgroup $C_2^3$ of $H$ is a
characteristic subgroup, hence $G$ acts on it via conjugation. This
defines a nontrivial homomorphism $f:G\to \Aut(C_2^3) \cong L_2(7).$ Let us treat
this case which also covers the third case for a possible quotient of
$G$.  Since $G = [G,G]$, the image of $f$ is not a solvable subgroup of
$L_2(7)$. The known classification of subgroups of $L_2(7)$ shows that
it is equal to $L_2(7)$. Thus the kernel  of $f$ is a subgroup of
order $2^3.3.5$. This order is in the List, but as we
will see in the next lemma it cannot be realized as the order of an
exceptional group. There is only one group with this order in Mukai's
list, which is  $\frakS_5$.
This group contains $\frakA_5$ as a unique
subgroup of index 2. The group $G$ acts on it by conjugation.
Since the group Out($\frakA_5)$ of outer automorphisms (modulo inner automorphisms) of
$\frakA_5$  is abelian, we get a nontrivial homomorphism $G\to
\frakA_5={\rm Inn}(\frakA_5)$. As
before we infer that it is surjective.  But we have seen in above that this is
impossible. This proves that $G$ is simple. It is known that there are two simple groups of order $2^6.3^2.5.7$, one is
$M_{21}$ and another is $\frakA_8$. The latter group contains an element of order 15 and must be excluded. The List gives only one possibility $\calR_G=A_2+2A_4+A_6+D_4$.

\medskip\noindent
Case 2. $G \ne [G,G]$.

\medskip
Since 5 and 7 divide $\#G$, $X/G$ has
singularities of type $A_4$ and $A_6$. This  follows from Table 1.
Assume $G$ has a normal subgroup $H$ with an abelian
quotient of order  5 (resp. 7).  Then $X/H \to X/G$ is a cyclic cover
of degree 5 (resp. 7). The pre-image of a singularity of type $A_6$
(resp. $A_4$) consists of 5 singularities of type $A_6$ (or 7
singularities of type $A_4$). This gives $\rank \calR_H > 20$, a contradiction.
This implies  that $35|\#[G,G]$. If $\mu([G,G]) \ge 5$, this contradicts
Proposition (4.2) from \cite{Mukai}. So $\mu([G,G]) \le 4$. On the other
hand, by Corollary (3.5) from \cite{Mukai},  $\mu([G,G])\ge \mu(G)= 4$.
Thus  $\mu([G,G])= 4$.  Replacing  $G$ with $[G,G]$ and repeating
the argument, we see that $G$ is obtained by taking extensions of
 a proper subgroup $K$ with $35|\#K$,
 $K=[K,K]$ and $\mu(K)=4$. By the result of Case 1, $K\cong M_{21}$ or
 $\frakA_7$. The first case can be excluded, as $M_{21}$ has order
 maximal in the List. The second case can also be excluded, as
 $\frakA_7$, though its order is not maximal in the List, admits no extensions
 in the set of finite symplectic groups. This easily follows from that
 $X/ \frakA_7$ has only one singularity of type $A_4$, and $5k,\ k>1,$
does not divide the order of a stabilizer subgroup.
\end{proof}

It follows from the List that a 2-group is not exceptional and hence
is in Mukai's list. We need the following description of
symplectic groups of order $2^6$.

\begin{lemma}\label{2^6} \cite{Mukai},\ \cite{Xiao}
There are 5 symplectic groups of order $2^6$, $\Gamma_{13}a_1$, $\Gamma_{22}a_1$,
$\Gamma_{25}a_1$, $\Gamma_{23}a_2$, and  $\Gamma_{26}a_2$.
\begin{enumerate}
\item[(i)]   $\Gamma_{23}a_2$ and  $\Gamma_{26}a_2$ do not contain a
  subgroup $\cong 2^4$.
\item[(ii)] $\Gamma_{13}a_1\cong 2^4:2^2\cong 4^2:2^2\cong 2^2.\ 2^4$.
\item[(iii)] $\Gamma_{22}a_1\cong 2^4:4\cong 2^3.\ (2\times 4)$.
\item[(iv)]  $\Gamma_{25}a_1\cong 2^4:2^2\cong 2^3.\ 2^3$.
\item[(v)] $\Gamma_{22}a_1$ and $\Gamma_{25}a_1$ have only one normal
  subgroup $\cong 2^4$.
\item[(vi)] $\Gamma_{13}a_1$,  $\Gamma_{22}a_1$ and $\Gamma_{25}a_1$
  always split over every normal subgroup $\cong 2^4$.
\end{enumerate}
In $(ii)-(iv)$, the last isomorphism is given by (commutator).(quotient).
\end{lemma}

\begin{proof} We give a proof of the last assertion (vi), which is not
  given explicitly in \cite{Mukai} or \cite{Xiao}.
By (ii)-(iv), each of these three groups has a  normal subgroup $\cong
2^4$
over which it splits. By (v) the assertion follows for
$\Gamma_{22}a_1$ and $\Gamma_{25}a_1$.

Let $G=\Gamma_{13}a_1$, and let $H_1,\ H_2$ be two distinct normal subgroups $\cong
2^4$ of $G$. If $\#(H_1\cap H_2)=8$, then the join
$H_1H_2$ is of order $2^5$. There is only one symplectic group of this
order containing a subgroup $\cong 2^4$ \cite{Mukai}. It is $\Gamma_{4}a_1$, but
this group
does not contain two subgroups $\cong 2^4$. This proves that $\#(H_1\cap
H_2)\le 4$. In this case $H_i$ contains a subgroup
$K_i\cong 2^2$ with $K_i\cap H_j=\{1\}$  ($\{i,j\}=\{1,2\}$). Thus
$H_j<G$ splits for $j=1,2$.
\end{proof}

\begin{corollary}\label{split}  Let $G$ be an exceptional group or in Mukai's
  list. Let $H$ be a normal subgroup of $G$. Assume either $H\cong
  2^4$ or $\#H=2^6$. Then $G$ splits over $H$.
\end{corollary}

\begin{proof}
Let $P$ be a 2-Sylow subgroup of $G$. Then $P$ contains $H$ as a
normal subgroup. By Gasch\"utz's theorem cited in \cite{Mukai} (p.204 in the proof of Proposition
(4.8)), it suffices to show that $H<P$ splits. We may assume $P\ne H$.

Assume $H\cong 2^4$. If $\#P=2^5$, then $P\cong\Gamma_{4}a_1$ and
$H<P$ splits. If $\#P=2^6$, then by Lemma \ref{2^6}, $H<P$ splits.
If $\#P=2^7$, then $P\cong F_{128}$, a unique symplectic group of
order 128. It follows from \cite{Mukai}
(Proposition (3.18)) that $H< P$ splits.

Assume $\#H= 2^6$. Then $P\cong F_{128}$, which splits over every
subgroup of index 2.
\end{proof}

\begin{lemma} Let $G$ be an exceptional group of order $2^a.3^b.5$ or
  $2^a.3.7$ with $a,b\ge 1$. Then $G$ is isomorphic to one of the
  following seven groups:
$$(2^3:7):3, \ \frakS_6,\ M_{10},\ 2^4:\frakA_6,\ 2^4:\frakA_5\ncong
 M_{20},\ M_{20}:2=2^4:\frakS_5,\ 2^3:L_2(7).$$
\end{lemma}

\begin{proof} Assume that $G$ is solvable. A solvable group of order
  $mn$, where $(m,n) = 1$ contains a subgroup of order
$m$ (see \cite{Hall}, Theorem 9.3.1).
If $\#G=2^a.3^b.5$, $G$ contains a subgroup of order $3^b.5$. Such
an order cannot be found in Mukai's list or in our List. The case
$\#G=2^6.3.7$ can also be excluded, because there is no group of order
$2^6.7$ in both lists. Assume $\#G=2^3.3.7$. If $G$ is simple, then
$G\cong L_2(7)$ with $\mu(L_2(7)) = 5$, hence, not exceptional. Thus
$G$ is not simple. Let
$H$ be a nontrivial normal subgroup of $G$. By inspecting possible
orders of $H$ and $G/H$, and using that an order 7 element does not
normalize an order 3 element, we infer that $\#H=2^3$ or $2^3.7$. In the
first case, the quotient $G/H$ has a normal subgroup of order 7, hence
we may assume the second case. A group of order $2^3.7$ does not
appear in Mukai's list, so must be exceptional. Later, we will see that
such a group must be isomorphic to $2^3:7$ and its lattice is
$3A_6+2A_1$. This implies that $G\cong (2^3:7):3$ and
$\calR_G=A_6+2A_5+2A_2$.
This order with the type of lattice can be found in the List.

Now we assume that $G$ is not solvable. Therefore $G$ contains a normal subgroup $H$
and a normal subgroup
$T$ of $H$ such that $H/T$ is a non-commutative simple group.

\medskip\noindent
Case 1. $T = \{1\}$, i.e. $G$ contains a simple nonabelian normal subgroup $H$.

\medskip
Assume  $5|\#G$. The known classification of simple groups of order
$2^a.3^b.5$ (see \cite{Brauer}) gives that $H$ is isomorphic to
$\mathfrak{A}_5$  or $\mathfrak{A}_6$.  Let $\bar{G} = G/H$. Since
the groups $\mathfrak{A}_5,\mathfrak{A}_6$ are in Mukai's list  and hence
satisfy $\mu(G) \ge 5$, we may assume that $\bar{G}$ is a nontrivial
group. It is known (see \cite{Xiao}) that  in both cases $X/H$ has
2 singular points of type $A_4$.
 Since $5k, k > 1,$ does not divide the order of a stabilizer subgroup of $G$,
we see that $\bar{G}$ acts simply transitively on the set of two
singularities, and hence is a cyclic group of order 2.
If $H \cong\mathfrak{A}_5$,
the group  $G$ is isomorphic to  $\frakS_5$ or the
direct product $\mathfrak{A}_5\times C_2$. The first group is in
Mukai's list and $\mu(G) = 5$. The second group contains a
cyclic group of order 10 which cannot act symplectically on $X$.

 So, we may assume that $H \cong \mathfrak{A}_6$. Again, $G$ cannot be
 the direct product $\frakA_6\times C_2$. The ATLAS \cite{ATLAS} shows that
$G\cong \frakS_6$ or $G \cong M_{10}$. The orders of these groups are
 in the List and there are 5 possible root lattices for groups of
 order $720$. Three of them contain sublattices of type $E_6$ or
 $D_5$. It is easy to see that the corresponding stabilizer subgroups
 isomorphic to $T_{24}$ or $Q_{12}$ are not subgroups neither of $S_6$
 nor $M_{10}$. This leaves only two possibilities $\calR_G =
 2A_3+A_4+2A_5$ or $\calR_G = A_2+A_3+A_4+A_7+D_4$ for groups of this
 order. Since $\frakS_6$ has no elements of order 8 but $M_{10}$ has,
 we see that the first case could be realized for $\frakS_6$ and the
 second for $M_{10}. $

 Assume $7|\#G$. The known classification of simple groups of order
$2^a.3.7$ \cite{Brauer} gives
that $H\cong L_2(7)$. This group is from Mukai's list and
$\mu(L_2(7)) = 5$.
So the quotient  $\bar{G} = G/H$ is a group of order
$2^{a-3}$. The orbit space $X/H$ has one singular point of
type $A_6$ which must be fixed under the action of $\bar{G}$ on $X/H$
giving a singularity $Gx$ on $X/G$ with stabilizer group $G_x$ of
order
$7.\#\bar{G}$. Table 1  shows that $\bar{G}$
must be trivial. So $\mu(G) = 5$ and this case is not realized.

 \medskip\noindent
 Case 2. $T\ne \{1\}$ and $H = G$.

\medskip
 Assume $5|\#G$. Again $\bar{G} = G/T\cong \frakA_5$ or $\frakA_6$. Thus $T$
is a normal subgroup  of order $2^{a-2}.3^{b-1}$ if
$\bar{G}\cong \frakA_5$, or $2^{a-3}$ if $\bar{G} \cong \frakA_6$.
Assume $b=2$ and $G/T = \frakA_5$.  Since an element $g_5$ of order 5
from $G$ does not normalize a subgroup of order 3, we see that the
number of Sylow 3-subgroups in $T$ must be divisible by 5. This is
a contradiction. Assume $T$ is of order $2^{a-2}$ and $G/T =
\frakA_5$.
Since an element of order 5 does not commute with an
element of order 2, it cannot act identically on  the center of $T$.
It is easy to check that an abelian group $A$ of order $2^n$ admits
an automorphism of order $5$ only if $n > 4$ or $A = C_2^4$. A group
which contains a central subgroup of index 2 is abelian. Thus $T$ is
either a 2-elementary abelian group of order $2^4$ (if $a = 6$), or
an abelian group of order $2^5$ (if $a = 7$).   It follows from Nikulin
(\cite{Mukai}, Proposition (3.20)) that an abelian group of order
$2^5$ does not admit the Mathieu representation. So $T \cong C_2^4$ . The
order $2^6.3.5$ is in our List.  There is only one possible root
lattice in this case $\calR_G = A_1+2A_4+D_5+E_6.$ This group is not
isomorphic to $M_{20}$, since $\mu(M_{20})=5$.  A possible scenario
is that $X/T$ has 15 singular points of type $A_1$ and $G/H \cong
\frakA_5$ has two orbits of 5 and 10 points on this set. One orbit
gives a singularity of type $E_6$, another one of type $D_5$.
By Corollary \ref{split}, $T<G$ splits, i.e. $G \cong
 2^4:\frakA_5$. There are two non-isomorphic actions of
$\frakA_5$ on $\bbF_2^4$ (see section 2.8 in \cite{Mukai}). One gives the
group $M_{20}$
from Mukai's list. Another one gives a group from (vii). We denote
this group by $M'_{20}$. It is realized as a subgroup of $M_{23}$.

 Assume $G/T = \frakA_6$. A similar argument shows that  $T \cong
 C_2^4$.  A direct computation shows that  $\mu(\frakA_6) = 5$. For
any nontrivial normal subgroup $T$ of $G$ we have $\mu(G/T) > \mu(G)$
 (\cite{Mukai}, Proposition (3.3)). Thus $\mu(G) \le 4$ and such group
 may appear. The order $2^4.\#\frakA_6$
 appears in the List with $\calR_G = A_2+A_3+2A_4+E_7$.  This case
 cannot be excluded. By Corollary \ref{split}, $T<G$ splits and $G$ is isomorphic to
 a semi-direct product $2^4:\frakA_6$.

Assume $7|\#G$. We have $\bar{G}\cong L_2(7)$ and $T$ is of order
$2^{a-3}$. It follows from the List that  $a = 6$ and
$\calR_G = A_2+A_3+A_5+A_6+D_4$ or $A_2+2A_3+A_6+E_6$.
Let us exclude the first possibility for $\calR_G$.  Since the
2-part of the discriminant group is generated by $> 2$ elements, by Lemma \ref{generators}, it
suffices to show that $G =[G,G]$.  Since $L_2(7)$ is simple, the image
of $[G,G]$ in $\bar{G}$ is  either  trivial or the whole $\bar{G}$. In
the first case $[G,G] \subset T$ and the abelian group $G/[G,G]$ maps
surjectively to a non-commutative group $L_2(7)$. In the second case,
$\#L_2(7)\le \#[G,G] \le 2^3\#L_2(7)$. The List shows that one of the
inequalities is the equality. If it is the first one, $G$ contains a
normal subgroup isomorphic to $L_2(7)$ and hence is isomorphic to the
product $T\times L_2(7)$ which contains elements of order larger
than allowed. If it is
the second one, then $G = [G,G]$. This proves that $\calR_G =
A_2+2A_3+A_6+E_6$. Next we prove that the extension $G=2^3.L_2(7)$
splits. From the type of $\calR_G $ we see that $G$ contains no elements of order
8. Let $P$ be a 2-Sylow subgroup of $G$. Since the 2-Sylow
subgroups of $L_2(7)$ are isomorphic to the dihedral group $D_8$, we
have the extension $P=2^3.D_8$. By Gasch\"utz's theorem cited in
\cite{Mukai} (p.204 in the proof of Proposition
(4.8)), it suffices to show that $2^3<P$ splits. Since $P$ is
isomorphic to one of
the five groups from Lemma \ref{2^6} and contains no elements of order
8, we see that $P$ must be isomorphic to
$\Gamma_{25}a_1$, $\Gamma_{23}a_2$, or  $\Gamma_{13}a_1$.
The last group can further be excluded as it does not have a normal
subgroup $\cong 2^3$ with quotient $\cong D_8$ \cite{HS}.
Consider a subgroup $K\cong 2^3.\frakS_4$ of $G$ containing $P$.
The Mukai's list contains 3 groups of the order $\#K$,
$T_{192}$, $H_{192}$,  and $\Gamma_{13}a_1:3$ \cite{Xiao}. The first two are
split extensions $2^3:\frakS_4$ (p.192 Remark \cite{Mukai}), and the last one cannot
contain $P$, as $P\ncong\Gamma_{13}a_1$. If $K$ is exceptional,
then as we will see later, $K$ is isomorphic to
(xi) $(Q_8\circ Q_8)^{\text{\Tiny \textbullet}}\frakS_3$,
(xxii) $(2^4:2)^{\text{\Tiny \textbullet}}\frakS_3$,
(xxiii) $\Gamma_{13}a_1:3$, or
(xxv) $\Gamma_{25}a_1:3$. The first two contain elements of order
8, as their $\calR_K$ show, and the third can also be excluded for the same
reason as above. Thus $P\cong\Gamma_{25}a_1$. Since $\Gamma_{25}a_1$
is a 2-Sylow subgroup of Hol$(2^3)=2^3:L_2(7)$  \cite{HS},
 $2^3<P$ splits.

\medskip\noindent
Case 3: $T\ne \{1\}, G\ne H$.

\medskip
It is shown in \cite{Mukai}, Theorem (4.9) that a non-solvable group
$G$ with $\mu(G)\ge 5$ is isomorphic to $\frakA_5, \frakA_6, \frakS_5,
L_2(7)$ or $M_{20} \cong 2^4:\frakA_5$. If $\mu(H) \ge 5$, we must
have $H \cong M_{20}$ and $T = C_2^4$. Since $X/M_{20}$ has two
singularities of type $A_4$, $G/H$ must be a cyclic group of order
2, and hence $G\cong M_{20}{}_{\text{\Tiny \textbullet}} 2$. The order $2\#M_{20}$ appears in the
List with $\calR_G= A_2+A_3+A_4+A_5+D_6$. From this lattice, it is
easy to compute the order breakdown for $G$, in particular $G$ has
more elements of order 2 than $M_{20}$. So the extension splits and
$G=M_{20}:2\cong 2^4:\frakS_5$.

Assume $\mu(H) = 4$. It follows from Case 2 that $H$ is isomorphic to
one of the two groups
$M'_{20}$, or $2^4:\frakA_6$. There are no orders in the
List of orders strictly divisible by the order of the latter group.
In the first case, as we saw
in Case 2, $X/M'_{20}$ has singularities of type $A_1+2A_4+D_5+E_6$. It has
2 singular points of type $A_4$, so $G/M'_{20}\cong
C_2$ and its action on $X/M'_{20}$ fixes the unique singularity of type
$D_5$. However a degree 2 extenson of $Q_{12}$ cannot be a stabilizer
subgroup. This proves that the group  $M'_{20}$ has
no extension in the set of finite symplectic groups.
\end{proof}

It remains to consider the cases where $\#G$ is divisible by at most
two primes. These are all solvable groups.

It follows from the List that
an exceptional group $G$ cannot be of order $q^a, q = 2,3,5,7$ or
$3^b5,\ 3^b7, b=1,2$.

\begin{lemma} Let $G$ be an exceptional group. Assume $G$ is solvable
of order $2^6\cdot 5$, or $2^3\cdot 7$, or $2^a\cdot 3^2, (3\le a\le7)$. Then $G$ is one
of the groups  from (ix)-(xvii) in Theorem \ref{except}.
\end{lemma}

\begin{proof}  First of all Proposition (5.1) from \cite{Mukai} gives
that a nilpotent $G$ is
either abelian with no elements of order 4 or a 2-group. By the order
condition
$G$ is not a 2-group. Assume that $G$ is abelian with no elements of
order 4. Then $G\cong C_2^a\times C_3^2, C_2^6\times C_5$, or
$C_2^3\times C_7$. The latter two cases have no Mathieu
representations. In the first case,
$\mu(G)=4$ implies that $a=2$. But such an order is not in the List.
 Thus we may assume that $G$ is not nilpotent and hence its Fitting
subgroup $F$ (maximal nilpotent normal subgroup) is a proper nontrivial
subgroup. We use Mukai's classification by analyzing only the cases
where the assumption $\mu(G) \ge 5$ was used.

The first  case is when  $F \cong C_3^2$.   The quotient $G/F$ must be
a 2-subgroup $H$ of the group $\Aut(C_3^2) \cong \GL_2(\bbF_3)$ of order
48.  The List shows that $\#H = 8$ or $16$. Assume that  $\#H = 8$.
There are 2 possible root
lattices for exceptional groups of order 72: $\calR_G =
A_1+A_2+A_3+2A_7$ or $D_4+2D_5+E_6$. It is known that $X/F$ has 8
singular points of type $A_2$. To get a singular point of type $E_6$
in $X/G$ the group $H$ must fix exactly one of the singular points
on $X/F$. This is obviously impossible. Thus only the first case is
realized. It corresponds to the case $G/F \cong C_8$ which leads to a
group with $\mu(G) = 4$ (see \cite{Mukai}, p. 206).  If $\# H = 16$,
the group $H$ is a 2-Sylow subgroup of $\GL_2(\bbF_3)$, known to be
isomorphic to the semi-dihedral group $SD_{16}$. The order of
$3^2:SD_{16}$ is 144 and it is in  the List with $\calR_G =
A_7+D_4+A_5+A_3+A_1$ or other 3 possibilities, all containing one copy
of the root lattice $E_6$ or $E_7$. These cases can
be easily excluded.  This gives
the groups from (ix) and (x).

The second possible case is when $G/F$ is of order divisible by 7 and $G$
has a group $K$ isomorphic to  $C_2^3:C_7$ as a subquotient group.
The only possible order from the List of the form $2^a.7$ is $2^3.7$.
The corresponding root lattice is $2A_1+3A_6$. This leads to our group in (xi).

The third case is when $\#G = 2^a.5$ and $G$ admits a quotient
$\bar{G}$ containing a subgroup $G_0$ with $\mu(G_0) = 5$ isomorphic
to $C_2^4: C_5$. The inspection of the List  gives that
$\#G = 2^6.5$. If $\bar{G} \ne G$, then $\mu(G) < \mu(G_0)$.
Assume this is the case. The kernel $H$ of the projection $G\to
\bar{G}$
is of order $2^s$ with $ s \le 2$. Since it is a normal subgroup
and an element of order 5 does not commute with an element of
order 2, the order of $\Aut(H)$ must be divisible
by $5$. This implies that $H =\{1\}$. Thus the Fitting
subgroup $F$ of $G$ is a normal subgroup isomorphic to $C_2^4$, and $G$ must contain
$C_2^4:C_5$ as a proper subgroup. It is
known that  $G/F$ is mapped injectively in $\textup{Out}(F)$.
 The quotient $G/F$ is of order
$2^2.5$. Since it cannot contain an element of order 10, the quotient is isomorphic to
$C_5: C_4$. There are 3 possible root lattices $\calR_G$ in the
List. Two of them contain a sublattice of type $D_4$ or
$D_6$. It is easy to see that the corresponding stabilizer subgroups
$Q_8, Q_{16}$ are not subgroups of $G$. The remaining case
$3A_3+A_4+A_7$ cannot be excluded. By Corollary \ref{split}, the extension $G=2^4.\ (5:4)$
splits, and gives case (xii).

The fourth case is when $G/F$ is of order divisible by 9. In this
case $G$ is of order  $2^a.3^2$ by the order assumption.
  Mukai considers  the Frattini subgroup $\Phi$ of $F$
(the intersection of maximal subgroups) and shows that
$F/\Phi \cong C_2^4$ and $G/\Phi$ contains a subgroup $G_0$ with
$G_0 \cong C_2^4:C_3^2 \cong \frakA_4\times \frakA_4$ and
$\mu(G_0) = 5$. As in the previous case, to get $\mu(G) = 4$  we
must have either $\Phi \ne \{1\}$ or $\Phi = \{1\}$ and $G_0$ is
a proper subgroup of $G$.

Assume $\Phi \ne \{1\}$ and let $\#\Phi= 2^s, s\le 3$. Assume
$s = 1$. The quotient $X/\Phi$ has 8 singular points of type $A_1$.
The group $G_0$ acts on $X/\Phi$ and permutes these points with
at least one stabilizer subgroup of order divisible by $3^2$. The known
structure of stabilizers shows that this is impossible. No stabilizer
$G_x$ is of order divisible by 9. Assume $s = 2$ thus $\Phi \cong
C_2^2$ or $C_4$.
In the first case the quotient $X/\Phi$ has 12 singular points of type
$A_1$ and the stabilizer
subgroups  of $G_0$ of these points are groups of order $2^2.3$ (one
orbit), or $2^3.3$ (two orbits of size 6 each). The first case gives
one singularity of type $E_6$, and the second 2 singularities of type
$E_7$. Both appear in the List. But, since $G_0$ has
singularities of type $2A_5+4A_2+A_1$ (no $A_3$)
on a minimal resolution $Y$ of $X/\Phi$ (see \cite{Xiao}) we easily
exclude the second case,  obtaining that the
singularities of $Y/G_0$ are of types $E_6+2A_2+2A_5$. This case can
be found in the List. This gives a possible case
$G \cong 2^2.\ (\frakA_4\times \frakA_4)\cong \Gamma_{13}a_1:3^2$ from (xiii).
If $\Phi = C_4$, $X/\Phi$ has singularities of type $4A_3+2A_1$.
The stabilizer of $G_0$ of a point of type $A_1$ must be a group
of order divisible by $2^3.3^2$.  No stabilizer $G_x$ is of order
divisible by 9. So this case does not occur.
A degree 2 extension of the group from (xiii) may also appear. In this
case $G/\Phi$ is a degree 2 extension of $G_0= \frakA_4\times \frakA_4$
and $\mu(G/\Phi)>\mu(G)=4$, so $G/\Phi$ appears on Mukai's list. There
appears only one such group in \cite{Xiao}, and it is $\frakA_{4,4}$. Thus
$G \cong 2^2.\ \frakA_{4,4}\cong \Gamma_{13}a_1:\frakA_{3,3}$, and we obtain that
$\calR_G=E_7+2D_5+A_2+A_1$. It occurs in the List.
This is the case (xiv).

Assume $s = 3$. Then $F$ is a Sylow 2-subgroup  of $G$  of order $2^7$
with normal subgroup $\Phi$ of order 8. Thus $\Phi\cong Q_8, D_8, 2^3,
2\times 4$ or $C_8$. Assume that
$\Phi \cong Q_8$. The quotient $X/\Phi$ has singular points
of type $2D_4+3A_3$ or $4D_4+A_1$. In any case the group
$G_0$ of order $2^4.3^2$ permutes points of type $D_4$ with
stabilizer of order divisible by 9. No stabilizers of this
order could occur. Similar argument rules out the remaining
possibilities $\Phi\cong D_8, 2^3, 2\times 4$ or $C_8$.
 This shows that the case $\#G = 2^7.3^2$ and $\#\Phi=8$ does not occur.

Assume $\Phi =\{1\}$. Thus $G$ contains $2^4: 3^2$ as
a proper subgroup and $F = 2^4$. If $a = 5$, the quotient group
$G/F$ is of order 18. There are two possible groups of order 18
which can act symplectically on $X$, a group isomorphic to
$\frakA_{3,3}$ or $C_3\times D_6$. The first case leads to a group
$G \cong \frakA_{4,4}$ with $\mu(G) = 5$. The second case gives
singularities of $X/G$ of types
$3A_5+A_3+A_2$. A case with a group of this order and the same $\calR_G$
can be found in the List. By Corollary \ref{split},
the extension $F<G$ splits, and gives the case (xv).
If $a = 6$, we have 3 possible groups $G/F$
isomorphic to $3^2:4$,
$3\times \frakA_4$, or $\frakS_{3,3}$. In the second case,
$G/F$ contains a normal subgroup of order 4, and its pullback is a
nilpotent normal subgroup of $G$ contaning $F$. This contradicts to
the maximality of the Fitting subgroup $F$. This excludes the second
group. The first and the third
groups have singularities of the quotient
of types $4A_3+2A_2+2A_1$ and $2A_5+A_2+6A_1$, resp..
The groups act on the set of 15 singular points of $X/F$.
It is easy to see that  in these two cases $X/G$ has singularities of types
$A_7+D_5+2A_3+A_2$, $ D_5+D_4+2A_5+A_1$. Both cases
can be found in the List. Both extensions split by Corollary
\ref{split}, giving the groups from (xvi) and (xvii).
If $a = 7$, we have 3 possible groups $G/F$
isomorphic to one of the following groups $\frakA_{4,3},\ N_{72} \cong 3^2.D_8,\
M_9 \cong 3^2.Q_8$ of order 72. The first case can be excluded by the
maximality of $F$, as the group contains a normal subgroup of order 4.
The order $2^7.3^2$ appears in the List with possible root lattices $\calR_G$ of types
$A_1+A_2+2D_5+E_7$ or $A_1+A_2+D_5+D_6+E_6$. It is easy to see that
the stabilizer subgroups $T_{24}$ and $O_{48}$ of singularities of
type $E_6$ and $E_7$ are not subgroups of $2^4_{\text{\Tiny\textbullet}} N_{72}$
or $2^4_{\text{\Tiny \textbullet}} M_9$. Thus these two groups are
also excluded.
\end{proof}

\begin{lemma}  Let $G$ be an exceptional group. Assume $G$ is solvable
of order $2^a\cdot 3$, $(4\le a\le7)$. Then $G$ is one
of the groups  from (xviii)-(xxvi) in Theorem \ref{except}.
\end{lemma}

\begin{proof} Again we follow the arguments from \cite{Mukai}.
It follows from the List that a 2-group is not exceptional and hence
is in Mukai's list. Hence his
assumptions (7.1) are satisfied except the last one where we have to
replace the condition $\mu(G) \ge 5$ with the condition $\mu(G)
\ge 4$.

First Mukai considers the case when the Fitting subgroup $F$ of $G$ is
of order divisible by 3, or equivalently, $G$ has a unique 3-Sylow subgroup
$T$. Let $S$ be a 2-Sylow subgroup and $\phi:S\to \Aut(T) = C_2$ be
the natural homomorphism. The classification of abelian nilpotent symplectic
groups shows that $\#\Ker(\phi)$ is of order $\le 4$, thus $\#G \le
24$. There are no exceptional groups of order $\le 24$.

Thus we may assume that the Fitting subgroup $F$ is a 2-group.
In this case $F$ is the intersection of all 2-Sylow subgroups of $G$.
If $F$ is equal to a unique 2-Sylow subgroup, then $G/F \cong
C_3$. Otherwise $G$ contains three 2-Sylow
subgroups and $G/F \cong \frakS_3$.

Using the classification of 2-groups, Mukai lists all possible groups
$F$. In our case they can be only of order $2^c,\ c=3,4,5,6,7$ ($c=a$, or
$a-1$).

 \smallskip\noindent
Case $c = 3$.

In this case $F \cong C_2^3$ or $Q_8$ and $G/F \cong \frakS_3$. There
 are 4 possible root lattices for exceptional groups of order 48. If $F
 \cong C_2^3$, the quotient $X/F$ has 14 singular points of type $A_1$
 so the largest stabilizer for the action of $G$ on $X$ is of order
 12. Since all possible root lattices contain a subdiagram $E_6$ or
 $E_7$ this case is not realized. The quotient $X/Q_8$ has 2 singular
 points of type $D_4$ and 3 singular points of type $A_3$, or other
 possibility is 4 points of type $D_4$ and one of type $A_1$. It is
 easy to see that the first possibility leads to the root lattice of
 type $A_2+D_4+2E_7$ and the second one to $A_2+D_5+D_6+E_7$. In both
 cases the group $G$ must coincide with the stabilizer of a singular
 point of type $E_7$ and hence is isomorphic
 to $O_{48}$. This is our case (xviii).

 \smallskip\noindent
Case $c = 4$.

 In this case $F \cong C_2^4, C_4^2$, or $Q_8\times C_2$, or
 $Q_8\circ C_4$.
In the last case, Mukai shows that $G\cong T_{24}\circ C_4$, or
$(T_{24}\circ C_4)_{\text{\Tiny \textbullet}} 2$, and excludes them because $T_{24}\circ C_4$
 contains an element of order 12.

Assume $F\cong Q_8\times C_2$.
This case is also excluded by Mukai
because of the assumption $\mu(G) \ge 5$. It cannot be excluded in
our case and leads to two groups $G =T_{24}\times C_2$, if $G/F \cong
C_3$, and  $G =(T_{24}\times C_2)_{\text{\Tiny \textbullet}} 2$, if $G/F \cong \frakS_3$. Let us
determine the root lattice of $G =T_{24}\times C_2$. The singular points of
$X/T_{24}$ are of type $E_6+D_4+A_5+2A_2$, or  $2E_6+A_3+2A_2$. In the
first case, the group $C_2$ fixes the unique singularity of type $E_6$.
Since $T_{24}\times C_2\neq O_{48}$, we get a contradiction. In the
second case $X/G$ has singularities of type $E_6+D_4+2A_5$. A group of
order $48$ with this root lattice is in the List and gives  case (xix).
It is easy to see that $G =(T_{24}\times C_2)_{\text{\Tiny \textbullet}} 2$ has root lattice
$2A_1+A_5+D_6+E_7$. Computing order breakdown, we see that $G$ has
more elements of order 2 than $T_{24}\times C_2$ or $O_{48}$. So the
extension splits and $G=(T_{24}\times C_2):2\cong O_{48}:2$. This is
also in the List  and gives  case (xx).

Assume $F\cong C_2^4$ or $C_4^2$.
If $F$ is a 2-Sylow subgroup, then $\#G = 48$. The quotient $G/F$ has either 15 singular
points of type $A_1$ or 6 singular points of type $A_3$. Thus the
largest possible order of a stabilizer subgroup of $G$ is $12$. The
List shows that the root lattice $\calR_G$ always contains a copy of
$E_6$ or $E_7$. This shows that this case does not occur.  So, $F$ is not a
2-Sylow subgroup and $G/F \cong \frakS_3$. There are 3 possible
root lattices for groups of order 96 in our List. They are of types
$2A_1+A_5+D_6+E_7, \ 2A_1+A_7+D_5+E_6,\ 3A_2+2A_7$. If $F \cong
C_2^4$, $X/F$ has 15 singular points of type $A_1$ and hence the
largest possible stabilizer subgroup of $G$ is of order 12 and no stabilizers
of order 8. This shows that this case is not realized. If $F \cong
C_4^2$,  $X/F$ has 6 singular points of type $A_3$ and hence the
largest possible order of a stabilizer subgroup of $G$ is 24. Since
$T_{24}\ncong 4_{\text{\Tiny \textbullet}}\frakS_3$, the root lattice of $G$ cannot contain
$E_6$. This rules out the first two root lattices, and the remaining root lattice is
$2A_7+3A_2$. The $3$-part of its discriminant group is isomorphic to
$3^3$. But $G\cong 4^2_{\text{\Tiny \textbullet}}\frakS_3$ does not admit a non-trivial
homomorphism to $C_3$. This
contradicts Lemma \ref{generators}. So this case must be excluded.

\smallskip\noindent
Case $c = 5$.

In this case $F$ is isomorphic to $Q_8\circ Q_8$ or
$2^4:2=\Gamma_4a_1$.

Assume  $F=Q_8\circ Q_8$. In this case $X/F$ has 9 singular points of type $A_1$ and 2 singular
points of type $D_4$. If $G/F \cong C_3$, the group $C_3$ fixes
the points of type $D_4$ and define 2 points of type $E_6$ on the
quotient $X/G$. No root lattices with such sublattices are realized
for groups of order 96. Thus $G/F \cong \frakS_3$.
There are 5 possible root lattices for
groups of order 192.  If $\frakS_3$ leaves points of type $D_4$
invariant, then $X/G$ has 2 singular points of type $E_7$. If it
permutes these point, we have one singular point of type $E_6$. By
inspection of the diagrams, we see only one lattice $2E_7+A_3+A_2+A_1$
 fits. From the lattice, it is easy to see that $G$ has exactly 19
 elements of order 2. Since the subgroup $Q_8\circ Q_8$ has the same
 number of elements of order 2, the extension does not split.
This is the group from (xxi).

Assume that $F \cong 2^4:2$. Then $X/F$ has 3 singular points
of type $A_3$ and 8 singular points of type $A_1$. Using the list of
possible root lattices for exceptional groups of order 96, we see that
$G/F \cong \frakS_3$. Since $T_{24}\ncong 4_{\text{\Tiny \textbullet}}\frakS_3$,
the group $\frakS_3$ does not leave a point of type $A_3$
invariant, and hence permutes  3 singular points
of type $A_3$, giving a singular point on $X/G$ of type $D_4$ or
$A_7$. Moreover, this argument shows that $X/G$ does not have a
singular point of type $E_7$ or $E_6$. By inspecting the
5 possible root lattices for groups of order 192,
we see only two lattices $A_1+A_2+A_7+2D_5$ and $A_1+3A_5+D_4$
survive. But the latter can also be eliminated by considering the
orbit decomposition of $\frakS_3$ on the 8 singular points of type
$A_1$.
This gives $\calR_G = A_1+A_2+A_7+2D_5$. By a similar computation of
the number of elements of order 2, we see that the extension does not
split. This is case (xxii).

\smallskip\noindent
Case $c = 6$.

In this case $F$ is a 2-group of type $\Gamma_{13}a_1$ or
$\Gamma_{25}a_1$  from Proposition (6.12) of \cite{Mukai}.

Assume $F \cong \Gamma_{13}a_1$. The quotient $X/F$ has 3
singular points of type $D_4$ and 6 singular points of type
$A_1$. If  $G/F \cong C_3$,  this leads to the two possible root lattices $A_1+D_4+3A_5$ and
 $2A_1+3E_6$. This is case (xxiii).
If $G/F \cong \frakS_3$, then we have a group of order
$2^7.3$. The corresponding root lattices in our List are of the following types
$A_1+A_3+A_5+D_5+D_6, A_1+A_3+A_5+D_4+E_7, A_1+2A_3+E_6+E_7$. The
first and the third correspond to the two root lattices of the
index 2 subgroup. By Corollary \ref{split}, the extension
$G=F.\ \frakS_3$ splits, and gives case (xxiv) and (xxiv$'$).

Assume $F \cong \Gamma_{25}a_1$. Then $X/F$ has  one singular point of
type $D_4$, three singular points of type $A_3$ and five singular
points of type $A_1$. If $G/F \cong C_3$, $\calR_G $ must
contains only one copy of  $E_6$. There is only one such lattice
$A_1+A_3+2A_5+E_6$ for groups of order 192.  This gives case (xxv).
If $G/F \cong \frakS_3$, the root lattice
$A_1+A_3+A_5+D_4+E_7$ may occur. Again, by Corollary \ref{split}, the
extension splits, and gives case (xxvi).

 \smallskip\noindent
Case $c = 7$.

In this case $G/F \cong C_3$ and Mukai leads this to contradiction.
It is still true in our situation. The group $F\cong F_{128}$ of order $2^7$
has singularities on $X/F_{128}$ of types $D_6+D_4+2A_3+3A_1$. Since $F$
is normal, the group $C_3$ acts  on $X/F_{128}$. It must fix the singular
points of type $A_3$ and gives on $X/G$ singular points of type $D_5$.
But $Q_{12}\ncong C_4.C_3$.
\end{proof}

\begin{remark}\label{subtle} Here we explain the difference between the groups
    $2^4:\frakS_{4}$.
Mukai's list contains a unique group of
order $384= 2^7.3$ which he denotes by $F_{384}$. Our list contains
three groups of this order non-isomorphic to $F_{384}$. The difference
between these groups and $F_{384}$ is very subtle. The group
$F_{384}$ is isomorphic to $2^4:\frakS_4$
and the pre-image of the normal subgroup $2^2$ of $\frakS_4$ is
isomorphic to $\Gamma_{13}a_1$ (see \cite{Mukai}, p.212). Thus
$F_{384}\cong \Gamma_{13}a_1:\frakS_3$ is an extension of the same
type as our groups (xxiv) and (xxiv$'$). The difference is of course in
the action of $\frakS_4$ on $2^4\setminus 0$, whose orbit
decomposition can be read off from the corresponding root lattice
$\calR_G$ and is given in the first column of Table 3.

\begin{table}[h]
\begin{center}
\begin{tabular}{|l||r|r|r|}
\hline
$G$&\small{action of $\frakS_4$ on $2^4\setminus \{0\}$}& \small{action of $\frakS_4$ on $\Xi$}&\small{composition series}\\ \hline
$F_{384}$& $3+12$&$1+1+2+4$&$(2,2)$\\ \hline
(xxiv)& $3+4+8$&$1+3+4$&$(2,1,1)$\\ \hline
(xxiv$'$)& $1+2+12$&$1+3+4$&$(1,1,2)$\\ \hline
(xxvi)& $1+6+8$&$1+1+6$&$(1,2,1)$\\ \hline
\end{tabular}
\end{center}
\caption{}
\end{table}

 Recall that in the action of $M_{24}$
on a set $\Omega$ of 24 elements,
the complement of an octad $\Xi$ is identified with the affine space $\bbF_2^4$,
and each element $g\in \frakA_8$ embeds in $M_{24}$ by acting on $\Xi$
as an even permutation and acting on the complement $\Omega\setminus \Xi$
as a linear map $i(g)$, where $i:\frakA_8\to L_4(2)=\GL_4(\bbF_2)$ is the exceptional
isomorphism of simple groups. For $F_{384}$ Mukai shows that the image
of $\frakS_4$ in $\frakA_8$ is a subgroup with orbit
decomposition $1+1+2+4$ (\cite{Mukai}, p.218 and Corollary (3.17)), dependent on the
assumption $\mu(G)\ge 5$. In our case, $\mu(G)=4$ and hence the number
of orbits must be equal to 3. It is easy to see that $1+3+4$ and $1+1+6$ are
the only possible
such decompositions for a subgroup of $\frakA_8$ isomorphic to $\frakS_4$.
On the other hand, it is known that $\frakS_4$ has only one non-trivial irreducible linear representation over a field of characteristic 2. It is isomorphic to the representation $\frakS_4\to \frakS_3 \cong L_2(2)$. This easily implies that any faithful linear representation of $\frakS_4\to L_4(2)$ is either decomposable as a sum of the trivial representation and an indecomposable  3-dimensional representation, or is an indecomposable representation with composition series with factors  of dimensions $(1,1,2), (2,1,1), (1,2,1), (2,2).$

In the reducible case, the representation is the direct sum of one-dimensional representation and an indecomposable 3-dimensional representation with composition series of type $(1,2)$, or $(2,1)$. In the first case $\frakS_4$ has 3 fixed points in $2^4\setminus \{0\}$. Assume that  the corresponding extension $G = 2^4:\frakS_4$ realizes. Then in the quotient $X/2^4$ we have 15 ordinary double points permuted by $\frakS_4$ according to its action on $2^4\setminus \{0\}$. Three of them are fixed. This implies that $X/G$ has 3 singular points of type $E_7$; this is too many. In the second case $\frakS_4$ has orbit decomposition of type  $1+3+3+4+4$. Again it is easy to see that this leads to a contradiction.
This argument can be also used to determine the second column of Table 3 using the known information about the root lattice $\calR_G$ (from \cite{Xiao} for a non-exceptional group or otherwise from the List.)

Thus we may assume that the representation of $\frakS_4$ in $\bbF_2^4$ is indecomposable. One can show that, up to isomorphism, the representation is determined by its composition series, hence we have 4 different cases. They correspond to the cases in Table 3, where the last column indicates the type of the composition series.

Note that the linear representations $\frakS_4\to L_4(2)$ from
(xxiv) and (xxiv$'$) differ by an exterior automorphism of
$L_4(2)$ defined by a correlation. The images of $\frakS_4$ under the
representations defined by  $F_{384}$ and (xxvi) are conjugate to
subgroups of $\frakS_6$ embedded in $\frakA_8 \cong L_4(2)$ as the
subgroup $\frakA_{2,4}$. They differ
by an exterior automorphism of $\frakS_6$.

The distinction between $M_{20}$ and  $M_{20}'$ also can be given
similarly; the image of $\frakA_5$ in $\frakA_8$ has orbit decomposition
$1+1+1+5$ and $1+1+6$, respectively.
\end{remark}

\begin{proposition}\label{subM23} All exceptional groups are contained in $M_{23}$.
\end{proposition}

\begin{proof}
We will indicate the chain of maximal subgroups starting from a maximal
subgroup of $M_{23}$ and ending at the subgroup containing the given group.
If no comments are given, the verification is straightforward using the
list of maximal subgroups of $M_{23}$ which can be obtained from ATLAS.
A useful fact is that in all extensions of type $2^4:K$, the group $K$
acts faithfully in $2^4$ (\cite{Mukai}, Proposition (3.16)) and hence
defines an injective homomorphism $K\to L_4(2)\cong\frakA_8$.

The subgroup $G(\Xi)$ of $M_{24}$ which preserves an octad $\Xi$ is isomorphic
to the affine group $\AGL_4(2) = 2^4:\frakA_8$. Here $2^4$ acts identically on the octad, and the quotient $\frakA_8$ acts on the octad by even permutations. A section $S$ of the semi-direct product is a stabilizer subgroup of a point outside the octad. Taking this point as the origin, we have an isomorpohism $i:\frakA_8\to S\cong L_4(2)$ which we used in Remark \ref{subtle}.  We will write any element from $G(\Xi)$
as $(g_1,g_2)$, where $g_1$ is a permutation of $\Xi$ and $g_2$ is a permutation of
$\Omega\setminus \Xi$. The group $\AGL_4(2)$ is generated by elements
$(g,i(g)), g\in \frakA_8,$ and translations $(1,t_a), a\in 2^4$.  The image of
$\frakA_8$ consists of elements $(g,i(g))$, where $i(g) \in L_4(2)$.  In particular,
all elements of $\frakA_8$ fix the origin in $\Omega\setminus \Xi$, and hence $\frakA_8$ is
isomorphic to a subgroup of $M_{23}$. Recall that the latter is
defined as the stabilizer of an element of $\Omega$.

(i) $M_{23}\supset \frakA_8$. We use that $\frakS_6$ embeds in
$\frakA_8$ as the subgroup $\frakA_{2,6}$.

(ii), (iii), (v), (vii) $M_{23} \supset M_{22}$.

(iv) $M_{23} \supset M_{21}\bull 2$.

(vi) $\subset$ (iii). We use the embedding of $\frakA_5$ in $\frakA_6$
as the  composition of the natural inclusion $\frakA_5\subset
\frakA_6$ and the exterior automorphism of $\frakS_6$.

(viii)  $M_{23}\supset M_{22}$. There is only one split extension $2^3:L_2(7)$.

(ix)  $\subset$ (x).

(x) $M_{23} \supset  M_{21}\bull 2 \supset (3^2:Q_8)\bull 2$.

(xi)  $\subset$ (xxvii).

(xii) $\subset$ (vii).
The group $5:4$ is the normailzer of a 5-Sylow subgroup of $\frakS_5$.
The  group $\frakS_5$ embeds in $\frakA_8 = L_4(2)$ as a subgroup $\frakA_{1,2,5}$
and its conjugacy class is unique. The conjugacy class of the subgroup $5:4$
is also unique, and  hence its action on $2^4$ is defined uniquely up to isomorphism.

(xiii)  $\subset$ (xiv).

(xiv) $M_{23} \supset 2^4:\frakA_7 \supset 2^4:\frakA_{3,4}= \Gamma_{13}a_1:\frakA_{3,3}$.

(xv) $\subset$  (xvii).

(xvii) $M_{23} \supset \frakA _8 \supset 2^4:\frakA_{2,3,3} \cong 2^4:\frakS_{3,3}$.

(xvi) $\subset$ (iii),  we use that $3^2:4$ is the normalizer of a 3-Sylow subgroup
of $\frakA_6$ and it is a unique (up to conjugation) maximal subgroup of $\frakA_6$.

(xviii)  $\subset$ (xxvi). First we use that
$O_{48}$ is a central extension of $\frakS_4$.  The group (xxvi) is an
extension  $2^4:\frakS_4$ and $\frakS_4$ fixes a unique point $a\in
2^4\setminus 0$ (see  Table 3). We embed $\frakS_4$ as a subgroup of $\L_4(2)$ generated by the matrices $x,y$ satisfying $x^4=y^2=(xy)^3=1$:
\[\begin{pmatrix}1&1&0&0\\
0&1&1&0\\
0&0&1&1\\
0&0&0&1\end{pmatrix},\quad
\begin{pmatrix}1&1&0&0\\
0&1&0&0\\
0&1&1&1\\
0&0&0&1\end{pmatrix}.\]
Let $a = (0,0,0,1), b = (1,1,1,1), c = (0,1,0,1)$. Then we immediately check that
$t_bi(x)=i(x)t_c,\, (t_bi(x))^4=(t_ci(y))^2=t_a$.
The subgroup generated by the  elements $(x,t_bi(x))$ and $(y,t_ci(y))$
is isomorphic to $O_{48}$.

(xviii), (xix) $\subset$ (xx). In particular, it gives another embedding of $O_{48}$ in $M_{23}$.

(xx)  $\subset$ (xxiv$'$).  The group (xxiv$'$) is an
extension  $2^4:\frakS_4$ and $\frakS_4$ fixes a unique point $a\in
2^4\setminus 0$ and leaves invariant a 2-dimensional subspace $<a,
a'>\subset 2^4$ (see  Table 3).
In this case, we realize $\frakS_4$ as a subgroup of $L_4(2)$ generated by the matrices
\[\begin{pmatrix}1&1&0&0\\
0&1&1&0\\
0&0&1&1\\
0&0&0&1\end{pmatrix},\quad
\begin{pmatrix}1&1&0&0\\
0&1&0&0\\
0&0&1&0\\
0&0&1&1\end{pmatrix}.\]
Let $a = (1,0,0,0), a' = (0,1,0,0), b = (1,1,1,1), c = (0,1,0,1) $. It is checked that $t_bi(x)=i(x)t_c,\, (t_bi(x))^4=(t_ci(y))^2=t_{a}$, and
the subgroup generated by the  elements $(x,t_bi(x))$, $(y,t_ci(y)),
(1,t_{a'})$ is isomorphic to $O_{48}:2$.

(xxi) $\subset$ (xxvi). Note $Q_8\circ Q_8=\Gamma_{5}a_1$, thus
the group (xxi) has a 2-Sylow subgroup $\cong
\Gamma_{26}a_2=\Gamma_{5}a_1\,^{\text{\Tiny \textbullet}}2$. If the
group (xxi) admits a degree 2 extension in our list of exceptional
groups or in the Mukai's list, it must be the group (xxvi). This
follows from the types of singularities. On the other hand, the group
(xxvi) contains a non-split extension of the from
$(Q_8\circ Q_8)^{\text{\Tiny \textbullet}}\frakS_3$. To show this,
denote  the groups (xxi) and (xxvi) by
$G$ and $K$, respectively.
The normal subgroup $\Gamma_{25}a_1$ of $K$ contains only one subgroup
$\cong \Gamma_{5}a_1$. Denote this subgroup by $H$. Then $H$ is normal
in $K$, and its quotient $K/H$ is of order 12 and contains a subgroup
$\cong\frakS_3$, hence $K/H\cong D_{12}$. Since
all symplectic groups of order $2^6$
are contained in the unique symplectic group of order $2^7$ which is
isomorphic to a 2-Sylow subgroup of $K$, there is
a chain of subgroups $H\subset\Gamma_{26}a_2\subset K$. This implies
that there is an order 2 subgroup $A$ of $K/H$ such that
$\phi^{-1}(A)$ is a non-split extension $\Gamma_{26}a_2=H^{\text{\Tiny
\textbullet}}2$, where $\phi:K\to K/H$ is the projection.
Since $K/H\cong D_{12}$, one can always find a
subgroup $B\subset K/H$ such that $A\subset B\cong\frakS_3$. Now
$\phi^{-1}(B)$ gives a non-split extension
$(Q_8\circ Q_8)^{\text{\Tiny \textbullet}}\frakS_3$.

(xxii)  $\subset$ (xiv).  Since $G=\Gamma_{4}a_1\,^{\text{\Tiny
    \textbullet}}\frakS_3$, it has a 2-Sylow subgroup $\cong
\Gamma_{22}a_1=\Gamma_{4}a_1\,^{\text{\Tiny \textbullet}}2=2^4:4$.
The normal subgroup $\Gamma_{4}a_1$ of $G$ contains only one subgroup
$\cong 2^4$. Denote this subgroup by $H$. Then $H$ is normal
in $G$, and its quotient $G/H$ is of order 12 and contains a cyclic subgroup
of order $4$, hence $G/H\cong Q_{12}$. Thus $G\cong 2^4:Q_{12}$, a
    split extension by Corollary \ref{split}. The group
$\frakA_{3,4}$ contains $<(12)(4567), (123)>\cong Q_{12}$.

(xxiii) $\subset$ (xxiv) or (xxiv$'$).

(xxiv), (xxiv$'$) $\subset$ (xiv).
Here and in the next inclusion use Table 3.

(xxv) $\subset$ (xxvi)  $\subset$ (iii).

(xxvii)  $\subset$ (viii). We use that
$L_3(2)\cong L_2(7)$ contains $7:3$ as the normalizer of a 7-Sylow subgroup.
\end{proof}


\begin{example}
The group $O_{48} \cong 2_{\text{\tiny \textbullet}}\frakS_4$ is contained in a maximal
subgroup of $U_3(5) = \PSU_3(\bbF_{25})$ isomorphic to $2_{\text{\tiny \textbullet}}\frakS_5$.
Thus it acts on a K3-surface from the example in the next section in the case $p = 5$, and its order is prime to $p$.
Unfortunately, we do not know how to realize explicitly other exceptional groups. It is known that the group $\frakA_5$ admits a symplectic action on the Kummer surface of the product of two supersingular curves in characteristic $p \equiv 2, 3 \mod 5$ \cite{Ibukiyama}. Together with the group $2^4$ defined by the translations, we have a symplectic group isomorphic to an extension $2^4:\frakA_5$. Unfortunately,  this group is $M_{20}$ not $M_{20}'$ as we naively hoped.
\end{example}

\section{Examples of finite groups of symplectic automorphisms
in positive characteristic $p$}

\medskip\noindent
{\it Non-exceptional tame groups of symplectic automorphisms}

\medskip
A glance at Mukai's list of examples of K3 surfaces with
maximal  finite symplectic group of automorphisms shows that all
of them can be realized over a field of positive characteristic $p
> 7$. A complete Mukai's list consists of 80 groups (81 topological
types) \cite{Xiao}. All of them realize over a field of positive
characteristic $p$, as long as the order is not divisible by $p$.

\medskip\noindent
{\it Wild groups of symplectic automorphisms}

\medskip
Here we give a list of examples of K3 surfaces over a field of
characteristic $p = 2,3,5,11$ with symplectic finite automorphism
groups of order divisible by $p$ not from the Mukai list.

$$ (p = 2)\quad  X = V(x^4+y^4+z^4+w^4+x^2y^2+x^2z^2+y^2z^2+xyz(x+y+z)) \subset \bbP^3$$
It  is a supersingular K3 surface with Artin invariant $\sigma =
1$ with symplectic action of the group PSL$(3, {\bbF}_4)_{\text{\Tiny \textbullet}} 2\cong
M_{21}{}_{\text{\Tiny \textbullet}} 2$, whose order is $(20,160). 2=(2^6. 3^2.
5. 7). 2$ (see \cite{DKo}). Although the order of this
group divides the order of $M_{23}$, it is not a subgroup of
$M_{23}$.

$$(p = 3) \quad X=V(x^4+y^4+z^4+w^4)\subset \bbP^3.$$ Fermat
quartic surface is supersingular
with Artin invariant $\sigma = 1$ in characteristic $p=3$ mod 4.
The general unitary group
GU$(4, {\bbF}_9)$ acts on the Hermitian form $x^4+y^4+z^4+w^4$
over ${\bbF}_9$, so PSU$(4,{\bbF}_9)=U_4(3)$ acts on $X$, which is
simple of order $3,265,920=2^7.3^6. 5. 7$. This
example was known to several people (A. Beauville,  S. Mukai,
T. Shioda, J. Tate). The order of this group does not
divide the order of $M_{23}$.

$$ (p = 5)\quad  X=V(x^6+y^6+z^6-w^2)\subset \bbP(1,1,1,3).$$
It is supersingular with Artin invariant $\sigma = 1$. The general
unitary group GU$(3, {\bbF}_{5^2})$ acts on the Hermitian form
$x^6+y^6+z^6$ over ${\bbF}_{5^2}$, so PSU$(3, {\bbF}_{5^2})=U_3(5)$, a simple group, acts
symplectically on $X$.  The order of this group is equal to $126,000=2^4.
3^2. 5^3. 7$ and does not divide the
order of $M_{23}$.

$$(p=7) \quad X=V(x^3+(y^8+z^8)x-w^2)\subset \bbP(4,1,1,6).$$
It is supersingular with Artin invariant $\sigma = 1$. The general
unitary group GU$(2, {\bbF}_{7^2})$ acts on the Hermitian form
$y^8+z^8$ over ${\bbF}_{7^2}$, so PSU$(2, {\bbF}_{7^2})\cong
L_2(7)$, a simple group of order 168, acts symplectically on $X$.
Although this group can be found in Mukai's list, the group action
on the surface in his example degenerates in characteristic 7.
This surface is birationally isomorphic to the affine surface
$y^2=x^3+(t^7-t)x$ (\cite{DK}, Examples 5.8). This surface is also
isomorphic to Fermat quartic surface in characteristic $p=7$, and
hence admits a symplectic action of the group $F_{384}$ of order 384.

$$(p=11) \quad X=V(x^3+y^{12}+z^{12}-w^2)\subset \bbP(4,1,1,6).$$
It is supersingular with Artin invariant $\sigma = 1$. The general
unitary group GU$(2, {\bbF}_{11^2})$ acts on the Hermitian form
$y^{12}+z^{12}$ over ${\bbF}_{11^2}$, so PSU$(2,
{\bbF}_{11^2})\cong L_2(11)$, a simple group of order $660$,
acts symplectically on $X$. This is a subgroup of $M_{23}$ but has 4
orbits, so it is not realized in characteristic 0. This surface is
birationally isomorphic to the affine surface $y^2=x^3+t^{11}-t$
(\cite{DK}, Examples 5.8).  This surface is also
isomorphic to Fermat quartic surface in characteristic $p=11$, and
hence admits a symplectic action of the group $F_{384}$ of order 384 (see Remark \ref{subtle}).

\bigskip

The last three examples were obtained through a discussion with S.
Kond\={o}. One can show  using an algorithm from \cite{Shioda2}
and its generalization \cite{Goto} that all the previous examples
are supersingular K3 surfaces with the Artin invariant equal to 1.


\end{document}